\documentclass[3p]{elsarticle}

\usepackage[nolist]{acronym}
\usepackage{hyperref}
\usepackage{amsfonts}
\usepackage{amssymb}  
\usepackage{graphicx,epstopdf}
\usepackage[caption=false]{subfig}
\usepackage{algorithmic}
\usepackage{acronym}
\usepackage{amsmath}

\journal{Journal of Control Engineering Practice}

\begin{acronym}
	\acro{UAV}{Unmanned Aerial Vehicle}
	\acro{CFD}{Computational Fluid Dynamics}
	\acro{AWE}{Airborne Wind Energy}
	\acro{OED}{Optimum Experimental Design}
	\acro{CAD}{Computer-Aided Design}
	\acro{IMU}{Inertial Measurement Unit}
	\acro{FCC}{Flight Control Computer}
	\acro{ODE}{Ordinary Differential Equation}
	\acro{NLP}{Non-Linear Program}
	\acro{CRLB}{Cramer-Rao Lower Bound}
	\acro{LTI}{Linear Time-Invariant}
	\acro{SNR}{Signal-to-Noise Ratio}
\end{acronym}

\newcommand{\VT}{V_\mathrm{T}}
\newcommand{\VTe}{V_{\mathrm{T}_\mathrm{e}}}

\begin{document}
\begin{frontmatter}
\title{Optimal Input Design for Autonomous Aircraft\tnoteref{acknowledgment}}
\tnotetext[acknowledgment]{This research was supported by the EU via ERC-HIGHWIND (259 166), ITN-TEMPO (607 957), ITN-AWESCO (642 682), by DFG via Research Unit FOR 2401, and by BMWi via the project eco4wind.}

\author[AP,UniFr]{G.~Licitra\corref{CorrespondingAuthor}}\ead{giovanni.licitra@imtek.uni-freiburg.de}
\author[UniKal]{A.~B\"{u}rger}\ead{adrian.buerger@hs-karlsruhe.de}
\author[AP]{P.~Williams}\ead{p.williams@ampyxpower.com}
\author[AP]{R.~Ruiterkamp}\ead{r.ruiterkamp@ampyxpower.com}
\author[UniFr]{M.~Diehl}\ead{moritz.diehl@imtek.uni-freiburg.de}

\cortext[CorrespondingAuthor]{Corresponding author}

\address[AP]{Ampyx Power B.V, The Hague, Netherlands}
\address[UniFr]{Department of Microsystems Engineering, University of Freiburg, Germany}
\address[UniKal]{Faculty of Management Science and Engineering, Karlsruhe University of Applied Sciences, Germany}

\begin{abstract}
Accurate mathematical models of aerodynamic properties play an important role in the aerospace field. In some cases, system parameters of an aircraft can be estimated reliably only via flight tests. In order to obtain meaningful experimental data, the aircraft dynamics need to be excited via suitable maneuvers. In this paper, optimal maneuvers are obtained for an autonomous aircraft by solving a time domain model-based optimum experimental design problem that aims to obtain more accurate parameter estimates while enforcing safety constraints.The optimized inputs are compared with respect to conventional maneuvers widely used in the aerospace field and tested within real experiments.
\end{abstract}
\begin{keyword}
Autonomous aircraft, optimum experimental design, optimization, system identification and parameter estimation.
\end{keyword}
\end{frontmatter}

\section{Introduction}

Nowadays, autonomous aircraft have become widespread for both civil and military applications. An important task for the development of these systems is mathematical modeling of the aircraft dynamics. Such models of aircraft dynamics regularly contain quantities called \textit{aerodynamic derivatives} (or simply \textit{derivatives}), which in general depend on the flight condition and the aircraft geometry.

The current practice is to retrieve derivatives from empirical data obtained from similar aircraft configurations or with tools based on \ac{CFD} and augmenting and verifying them by wind tunnel tests. For standard aircraft configurations such methods are generally in good agreement with experimentally obtained values. However, for less conventional configurations these tools provide only a rough approximation of the aerodynamic properties \cite{mulder2000flight}.

\begin{figure}[tbhp]
	\centering
	\includegraphics[width = 365pt, height = 120pt]{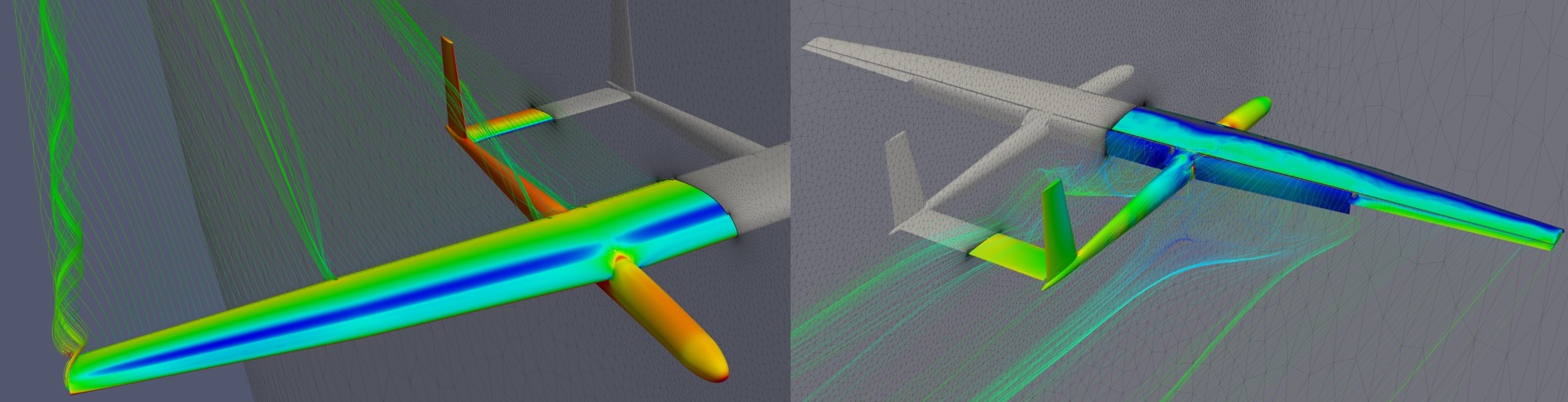}%
	\caption{CFD of the $3^{rd}$ prototype high lift, rigid wing autonomous aircraft designed by Ampyx Power B.V.}
	\label{fig:AP3CFD}
\end{figure}

This problem often arises in the \ac{AWE} community \cite{fagiano2012airborne,diehl2013airborne} where non-conventional high lift aircrafts need to be designed for extremely challenging operational environments \cite{Cap26AWEbook}. Figure~\ref{fig:AP3CFD} shows the \ac{CFD} analysis of a non-conventional high lift aircraft designed by Ampyx Power B.V. \cite{AP}. In this case, intensive flight test campaigns must be set in order to gain additional insight about the aerodynamic characteristics. 

A successful flight test campaign depends on many factors, such as selection of instrumentation, signal conditioning, flight test operations procedure, parameter estimation algorithm and signal input design. In this paper, the focus is on the optimization of signal inputs that aim to maximize the information content of the measurements data used for determining the aircraft aerodynamic properties via a model-based, time domain approach. A steady reference condition is considered and constraints are enforced in order to prevent flight envelope violation. The optimized experiments are assessed and a real flight test campaign is carried out. 

The paper is organized as follows. In section~\ref{sec:MathematicalModel}, the mathematical model of the rigid wing aircraft for \ac{OED} purposes is introduced and a brief overview of the case study designed by Ampyx Power B.V. is provided. Section~\ref{sec:FightTestProcedure} describes the flight test operation procedure and safety requirements as well as underlying theoretical and practical aspects. Section~\ref{sec:OEDproblem} presents the formulation of the \ac{OED} problem based on the Cramer-Rao Lower Bound. In section~\ref{sec:Results} the performance of the optimized inputs are assessed and compared to the performance of conventional input signals widely used in the aerospace field. The optimal maneuvers are first analyzed via a reliable flight simulator in section~\ref{sec:RealImplementation} and subsequently experimental data of a real flight test are provided. Section~\ref{sec:Conclusions} concludes.

\section{Mathematical Model}\label{sec:MathematicalModel}

In this section, the mathematical model of a rigid wing aircraft is introduced and a brief overview of the case study is provided.

\subsection{Modeling of Aircraft}
\begin{figure}[tbhp]
	\centering
	\includegraphics[scale = 0.33]{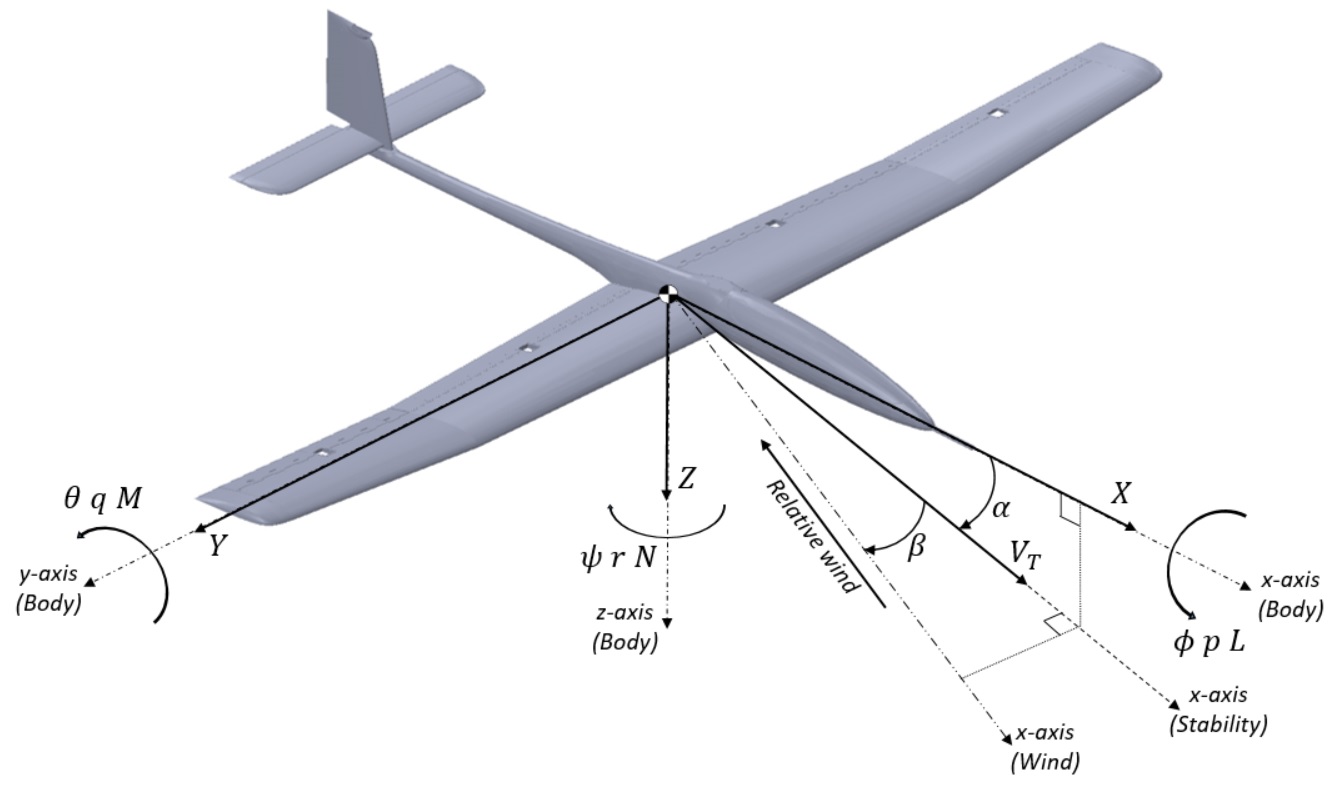}
	\caption{Definition of axes, Euler angles, aerodynamic states, forces and moments on a rigid wing aircraft.}
	\label{fig:AircraftConvention}
\end{figure}

For system identification purposes, let us consider the aircraft dynamics as follow \cite{stevens2015aircraft}
\begin{subequations}\label{eq:MathematicalModel}
\begin{align}
		\dot{V}_{\mathrm{T}}  & = \frac{Y\sin\beta + X\cos\alpha \cos\beta + Z \cos\beta  \sin\alpha}{m} + G_{\VT}, \label{eq:Vt}\\
		\dot{\beta}  & = \frac{Y\cos\beta  - X\cos\alpha \sin\beta - Z \sin\alpha \sin\beta}{m \VT} + \frac{G_{\beta}}{\VT} - r \cos \alpha + p \sin \alpha, \label{eq:beta}\\
		\dot{\alpha} & =  \frac{Z\cos\alpha - X\sin\alpha}{m \VT \cos\beta} + \frac{G_{\alpha}}{\VT \cos\beta} + \frac{q\cos\beta - (p\cos \alpha + r \sin \alpha)\sin\beta}{\cos\beta}, \label{eq:alpha}\\
		\dot{\phi}   & = p + r \cos\phi\tan\theta + q\sin\phi\tan\theta, \label{eq:phi}\\
		\dot{\theta} & = q\cos\phi - r\sin\phi, \label{eq:theta}\\
		\dot{\psi}   & = \sec\theta \left(q \sin\phi +  r\cos\phi \right), \label{eq:psi}\\
		\dot{p}      & = \frac{J_{xz}}{J_{x}} \dot{r} - qr\frac{\left(J_{z} - J_{y} \right)}{J_{x}} + qp \frac{J_{xz}}{J_{x}} + \frac{L}{J_{x}}, \label{eq:p}\\
		\dot{q}      & =  -pr \frac{J_{x} - J_{z}}{J_{y}} - (p^{2} - r^{2})\frac{J_{xz}}{J_{y}} + \frac{M}{J_{y}}, \label{eq:q}\\
		\dot{r}      & = \frac{J_{xz}}{J_{z}} \dot{p} -pq \frac{J_{y}-J_{x}}{J_{z}} - qr\frac{J_{xz}}{J_{z}} + \frac{N}{J_{z}}, \label{eq:r}	
\end{align}
\end{subequations}

where $(\VT,\beta,\alpha)$ are the aerodynamic states, i.\,e., true airspeed $\VT$, angle of side-slip $\beta$ and angle of attack $\alpha$, whereas the states $(\phi,\theta,\psi)$ denote the Euler angles of roll, pitch and yaw with $(p,q,r)$ the corresponding angular body rates \cite{stevens2015aircraft}.
The aircraft is assumed to have a constant mass $m$, moments of inertia $J_{x},J_{y},J_{z}$, cross moment of inertia $J_{xz}$ and to be subject to external aerodynamic forces $(X,Y,Z)$, moments $(L,M,N)$ and, obviously, gravity. More precisely, the gravity components are expressed as 
\begin{subequations}\label{eq:GravityComponents}
	\begin{align}	    
		G_{\VT}  & =  g_{D} \left(\sin\beta \sin\phi \sin\theta -\cos\alpha\cos\beta\sin\theta +  \sin\alpha \cos\beta \cos\phi \cos\theta \right), \\
		G_{\beta}  & =  g_{D} \left( \cos\alpha \sin\beta \sin\theta + \cos\beta\sin\phi \cos\theta - \sin\alpha \sin\beta \cos\phi \cos\theta  \right),\\
		G_{\alpha} & =  g_{D} \left( \sin\alpha \sin\theta + \cos\alpha \cos\phi \cos\theta  \right),
	\end{align}
\end{subequations}
with $g_{D} \approx 9.81\,\mathrm{m/s^{2}}$ the gravitational acceleration. The nomenclature introduced above is summarized in fig.~\ref{fig:AircraftConvention}. Furthermore, the mathematical model \ref{eq:MathematicalModel} implicitly presumes that the vehicle is a rigid body with a plane of symmetry such that the moments of inertia $J_{xy},J_{xz}$ are zero, whereas the Earth is assumed flat and non-rotating with a constant gravity field \cite{mulder2000flight}.

\subsection{Aerodynamic model}\label{sec:aerodynamics}

In flight dynamics, there are different methods of aerodynamic derivatives modeling. In many practical cases the aerodynamic forces and moments are approximated by linear terms in their Taylor series expansion. Such approximations yield sufficient accuracy for attached flows at low angles of attack \cite{etkin2012dynamics}.
In this case, the aerodynamic properties can be normalized with respect to the dynamic pressure $\bar{q} = \frac{1}{2} \rho \VT^{2}$ with $\rho\approx 1.225 \,\mathrm{kg/m^{3}}$ the free-stream mass density, and a characteristic area for the aircraft body
\begin{subequations}\label{eq:Fa_Ma}
	\begin{align}
		X & = \bar{q}S        \,C_{X} & Y & = \bar{q}S        \,C_{Y} & Z & = \bar{q}S        \,C_{Z} \\
		L & = \bar{q}Sb       \,C_{l} & M & = \bar{q}S\bar{c} \,C_{m} & N & = \bar{q}Sb       \,C_{n} 				
	\end{align}
\end{subequations}
In \ref{eq:Fa_Ma} $S,\,b,\,\bar{c}$ are reference wing area, wing span and mean aerodynamic chord, respectively, while $C_{X},\,C_{Y},\,C_{Z}$ denote the forces and $C_{l},\,C_{m},\,C_{n}$ the moment coefficients.
For conventional aircraft, the aerodynamic coefficients are usually broken down into a sum of terms as in \ref{eq:AdimensionaAerodynamicCoefficients}
\begin{subequations}\label{eq:AdimensionaAerodynamicCoefficients}
	\begin{align}
		C_{X} & = C_{X_{\alpha}}\alpha + C_{X_{q}}\hat{q} + C_{X_{\delta_{\mathrm{e}}}} \delta_{\mathrm{e}} + C_{X_{0}},\\
		C_{Y} & = C_{Y_{\beta}} \beta  + C_{Y_{p}}\hat{p} + C_{Y_{r}} \hat{r} + C_{Y_{\delta_{\mathrm{a}}}} \delta_{\mathrm{a}} + C_{Y_{\delta_{\mathrm{r}}}} \delta_{\mathrm{r}}, \\  
		C_{Z} & = C_{Z_{\alpha}}\alpha + C_{Z_{q}}\hat{q} + C_{Z_{\delta_{\mathrm{e}}}} \delta_{\mathrm{e}} + C_{Z_{0}},\\     	
		C_{l} & = C_{l_{\beta}} \beta  + C_{l_{p}}\hat{p} + C_{l_{r}} \hat{r} + C_{l_{\delta_{\mathrm{a}}}} \delta_{\mathrm{a}} + C_{l_{\delta_{\mathrm{r}}}} \delta_{\mathrm{r}}, \label{eq:Cl}\\
		C_{m} & = C_{m_{\alpha}}\alpha + C_{m_{q}}\hat{q} + C_{m_{\delta_{\mathrm{e}}}} \delta_{\mathrm{e}} + C_{m_{0}}, \\
		C_{n} & = C_{n_{\beta}} \beta  + C_{n_{p}}\hat{p} + C_{n_{r}} \hat{r} + C_{n_{\delta_{\mathrm{a}}}} \delta_{\mathrm{a}} + C_{n_{\delta_{\mathrm{r}}}} \delta_{\mathrm{r}},
	\end{align}
\end{subequations}
that depend on the normalized body rates $\hat{p}= \frac{b \, p}{2\VT},\hat{q} = \frac{\bar{c} \, q}{2\VT},\hat{r} = \frac{b \, r}{2\VT}$, angle of attack $\alpha$ and side slip $\beta$, as well as the control surface deflections which in this case are aileron $\delta_{\mathrm{a}}$, elevator $\delta_{\mathrm{e}}$ and rudder $\delta_{\mathrm{r}}$. The coefficients $C_{i_{j}}$ with $i = \{X,Y,Z,l,m,n\}$ and $j = \{\alpha,\beta,p,q,r,\delta_{\mathrm{a}},\delta_{\mathrm{e}},\delta_{\mathrm{r}},0\} $ are the \textit{dimensionless aerodynamic derivatives} that need to be identified. 
One has to highlight that these coefficients are typically valid only for small neighborhoods of a specific flight condition, hence system identification flight tests are usually performed for multiple flight configurations.

\subsection{Aircraft airframe}\label{sec:UAVairframe}

\begin{figure}[tbhp]
	\centering
	\includegraphics[scale = 0.6]{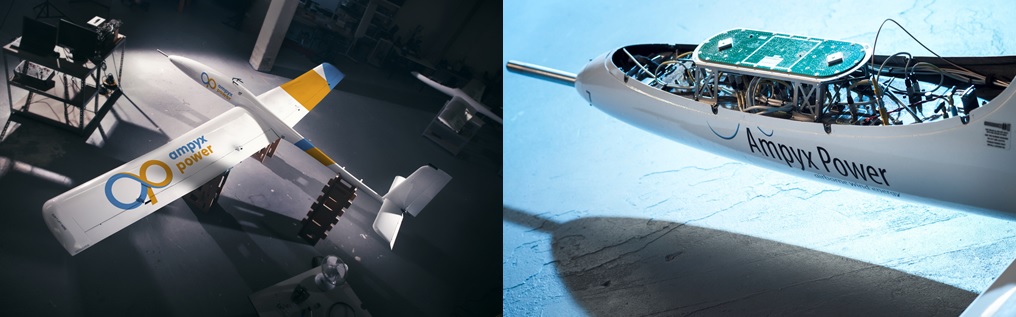}
	\caption{The $2^{nd}$ prototype high lift, rigid wing autonomous aircraft designed by Ampyx Power B.V. (left) and its on-board flight computer (right).}
	\label{fig:AP2}
\end{figure}

The case study considered within this paper and shown in fig.~\ref{fig:AP2} is a high lift, rigid wing autonomous aircraft used as airborne component of a pumping mode airborne wind energy system. Details on the system can be found in \cite{AWEbook,Cap26AWEbook}.
The control surfaces of the airframe are aileron $\delta_{\mathrm{a}}$, elevator $\delta_{\mathrm{e}}$ and rudder $\delta_{\mathrm{r}}$. Also, the airframe is equipped with flaps $\delta_{\mathrm{f}}$ which are not used within this analysis. All control surfaces are actuated via electric servos with $\delta_{\mathrm{a}} \in [-20\,^\circ, 20\,^\circ]$, $\delta_{\mathrm{e}}  \in [-30\,^\circ, 30\,^\circ]$ and $\delta_{\mathrm{r}} \in [-30\,^\circ, 30\,^\circ]$.
Further, the aircraft is equipped with a propulsion system which consists of two electric motors that drive two blades mounted on top of the fuselage. 
Physical properties of the airframe are summarized in table~\ref{tab:AP2_parameters} with moments of inertia obtained through \ac{CAD} models and validated via swing tests \cite{de1987accurate,lyons2002obtaining}.

The aircraft is instrumented with an \ac{IMU} which provides measurements of body angular rates and translational accelerations. The aerodynamic states are estimated by mean of a \textit{five hole pitot tube} which is mounted at the nose of the fuselage. The sensor noise of each component is expressed in terms of its standard deviation $\sigma_{y}$ as shown in table~\ref{tab:std_sensors}.
The actuator commands are delivered by an on-board flight computer at $100\,\mathrm{Hz}$, data are recorded at the same rate. The aircraft systems are battery powered, allowing $\approx 15\,\mathrm{min}$ of flying time on a single charge.

A flight test campaign that aims to identify the aerodynamic properties of an aircraft requires accurate models of the system dynamics that contains all contributions, e.\,g., thrust dynamics. Though, it is usually rather difficult to achieve high accuracy on propeller dynamics. Furthermore, the rotation of the blades introduces additional noise for each angular rate and acceleration channel. Therefore, in order to obtain meaningful experimental data for the aerodynamic parameter identification, the propulsion system is not taken into account in the mathematical model introduced in \ref{sec:MathematicalModel} and switched off whenever a signal excitation is conducted in order to discard the thrust effect from the aerodynamic forces and moments.

\begin{table}[tbhp]
	\caption{Physical properties of the aircraft designed by Ampyx Power B.V.}
	\label{tab:AP2_parameters}
	\centering
	\begin{tabular}{lccc}\hline
		\bf Name                &\bf Symbol &\bf Value & \bf Unit                  \\\hline
		mass                    & $ m $     & $36.8$   & $\mathrm{kg}$             \\
		moment of inertia       & $J_{x}$   & $25$     & $\mathrm{kg \cdot m^{2}}$ \\
		moment of inertia       & $J_{y}$   & $32$     & $\mathrm{kg \cdot m^{2}}$ \\
		moment of inertia       & $J_{z}$   & $56$     & $\mathrm{kg \cdot m^{2}}$ \\
		cross moment of inertia & $J_{xz}$  & $0.47$   & $\mathrm{kg \cdot m^{2}}$ \\
		reference wing area     & $S$       & $3$      & $\mathrm{m^{2}}$          \\
		reference wing span     & $b$       & $5.5$    & $\mathrm{m}$              \\
		reference chord         & $\bar{c}$ & $0.55$   & $\mathrm{m}$              \\\hline		
	\end{tabular}
\end{table}

\section{Flight test procedure and the rationale behind}\label{sec:FightTestProcedure}

This section describes the flight test procedure for autonomous aircraft augmented with safety, practical and theoretical aspects.

\subsection{Experimental setup and safety requirements}

The flight test procedure is divided into three parts: take-off, execution of the experiments and landing. The take-off and landing are performed manually by the pilot in remote. After the manual stabilization of the aircraft dynamics, the autonomous mode is enabled and the vehicle flies into a predefined race track pattern in order to establish the so called \textit{steady wing-level flight} condition (see section~\ref{sec:steadyConditionAndDecoupling}). From a safety point of view, typical flight tests for \ac{UAV} are limited to \textit{line-of-sight} range in order to both avoid communication dropout and to guarantee that the pilot can regain manual control of the aircraft at any time \cite{dorobantu2013system}.

To prevent biases due to correlation between the measurement noise and the inputs, it is best to perform open-loop experiments \cite{soderstrom2001identification} which are also beneficial for systems equipped with sensors that are susceptible to high levels of measurement noise \cite{dorobantu2013system}, in our case the five hole pitot tube which provide the aerodynamic states ($\VTe,\beta,\alpha$).
Immediately after the turn in the race track pattern is completed and the steady wing-level flight is achieved, the signal excitation is performed along one single axis. All remaining control surfaces are fixed at their trim values throughout the experiment.
Experiments are generally repeated on each axis so as to both obtain a rich data set and reduce the effect of sensor biases as well as colored noise (atmospheric turbulence) on estimation results \cite{licitra2017pe}.

Due to limit of the flight test field and steady state requirements prior to and after each turn, a $10\,\mathrm{s}$ experiment time window is chosen for the system identification flight test. As a consequence, dynamics at frequencies below $0.1\,\mathrm{Hz}$ cannot be identified accurately.

Finally, one has to point out that the autonomous mode of the aircraft implies no action of the pilot during the flight test unless system failures are detected. As a result, before conducting any real experiment, excitation signals must be tested within reliable simulations to reduce the chance of flight envelope violation and loss of track. Nonetheless, it may happen that during the real flight test the aircraft violates the flight envelope, e.\,g., due to significant inaccuracies of the a priori models or unexpected gust occurring during the open loop-phase. For this reason, \textit{flight envelope limit detection} algorithms should be programmed in the \ac{FCC} in order to avoid damaging or even destruction of the vehicle. 

\subsection{Steady wing-level flight condition and decoupling of dynamics}\label{sec:steadyConditionAndDecoupling}
From a both physical and practical point of view steady wing level flight condition is crucial for system identification purposes. More specifically, an aircraft is meant in steady wing-level flight when its body angular rate ($p,q,r$) and roll angle $\phi$ are equal to zero and it flies with constant airspeed $\VTe$ \cite{stevens2015aircraft}. Fulfillment of this steady condition allows decoupling of the aircraft motion in \textit{longitudinal} and \textit{lateral} dynamics, hence one can focus only on a subset of the entire aircraft dynamics which is mainly excited from a given maneuver.

For aircraft parameter estimation experiments, typically a linear perturbation model structure is assumed \cite{klein1990optimal}. Consequently, the flight test inputs need to be perturbations with respect to the steady condition so that the system response can be adequately modeled by such linear structure. Behind historical reasons for deriving the small-perturbation equations, one can obtain great insight into the relative importance of the various aerodynamic derivatives under different flight conditions and their effect on the aircraft stability \cite{stevens2015aircraft}.  

The longitudinal dynamics is described via LTI state-space form by the state $\mathbf{x_{lon}} = \left[V_{T}~\alpha~\theta~q\right]^\top$, which corresponds to \ref{eq:Vt}, \ref{eq:alpha}, \ref{eq:theta}, \ref{eq:q}. The forces $X$, $Z$ and the moment $M$ are assumed to be linear functions of $\VT,\alpha,q$ and the elevator deflection $\delta_{\mathrm{e}}$, resulting in the following system 
\begin{equation}\label{eq:longDyn}
\begin{bmatrix}
\dot{V}_{T}   \\
\dot{\alpha}  \\
\dot{\theta}  \\
\dot{q}  
\end{bmatrix} = 
\begin{bmatrix}
X_{V} & X_{\alpha}                   & -g_{D}\cos\theta_{e} & X_{q} \\
Z_{V} & \frac{Z_{\alpha}}{V_{T_{e}}} & -g_{D}\sin\theta_{e} & Z_{q} \\
0   &     0                        &           0          &  1    \\
M_{V} & M_{\alpha}                   &           0          & M_{q}
\end{bmatrix}
\begin{bmatrix}
V_{T}  \\
\alpha  \\
\theta  \\
q  
\end{bmatrix} + 
\begin{bmatrix}
X_{\delta_{e}}  \\
\frac{Z_{\delta_{e}}}{V_{T_{e}}} \\
0               \\
M_{\delta_{e}}  
\end{bmatrix} 
\delta_{e}
\end{equation}
where the non-zero elements are known as \textit{dimensional aerodynamic derivatives} while $\theta_{e}$ is the steady-state pitch angle. The dimensional derivatives can be converted in dimensionless derivatives shown in \ref{eq:AdimensionaAerodynamicCoefficients} via the geometrical configuration of the aircraft, for details see \cite{stevens2015aircraft,mulder2000flight}. 
The longitudinal dynamics can be further decoupled into the \textit{Phugoid} and \textit{Short-period} mode. The Phugoid mode is normally rather slow, slightly damped, and dominates the response in $\VT$ and $\theta$ while the Short-period mode is typically fast, moderately damped, and dominates the response in $\alpha$ and $q$. For control applications, accurate knowledge of the Phugoid mode is not crucial due to the low frequency of the oscillation which is compensated via feedback control, whereas the Short-period mode is crucial for stability and performance characteristics \cite{cook2007flight}. 

The lateral dynamics are described analogously by the state $\mathbf{x_{lat}} = \left[\beta~\phi~\psi~p~r\right]^\top$, which corresponds to equations \ref{eq:beta}, \ref{eq:phi}, \ref{eq:p}, \ref{eq:r}. Force $Y$ and moments $L$ and $N$ are described by linear functions of $\beta,p,r$ and inputs $\mathbf{u_{lat}} = \left[ \delta_{a}~\delta_{r}\right]^\top$. The resulting system is given by
\begin{equation}\label{eq:latDyn}
\begin{bmatrix}
\dot{\beta} \\
\dot{\phi}  \\
\dot{p}     \\
\dot{r}  
\end{bmatrix} = 
\begin{bmatrix}
\frac{Y_{\beta}}{V_{T_{e}}}  & g_{D}\cos\theta_{e} & Y_{p}  & Y_{r} - V_{T_{e}} \\
0              &      0              &   1    & \tan \theta_{e}   \\
L_{\beta}'         &      0              & L_{p}' & L_{r}'            \\
N_{\beta}'         &      0              & N_{p}' & N_{r}' 
\end{bmatrix}
\begin{bmatrix}
\beta \\
\phi  \\
p     \\
r   
\end{bmatrix} + 
\begin{bmatrix}
\frac{Y_{\delta_{a}}}{V_{T_{e}}}  & \frac{Y_{\delta_{r}}}{V_{T_{e}}} \\
0             &                     0            \\
L_{\delta_{a}}'      &              L_{\delta_{r}}'     \\
N_{\delta_{a}}'      &              N_{\delta_{r}}'
\end{bmatrix} 
\begin{bmatrix}
\delta_{a} \\ 
\delta_{r}         
\end{bmatrix} 
\end{equation}
where $\psi$ has been dropped and its derivatives are defined in \cite{mcruer2014aircraft}. Unlike the longitudinal dynamics, the lateral dynamics cannot be decoupled into independent modes. They are governed by a slow \textit{Spiral} mode, a fast lightly damped \textit{Dutch roll} mode, and an even faster \textit{Roll Subsidence} mode. The Spiral and Roll Subsidence mode usually involve almost no side-slip $\beta$. The Roll Subsidence mode is almost pure rolling motion around the x-axis whereas the Spiral mode consist of yawing motion with some roll. Often, modern aircraft have an unstable Spiral mode which cause an increment of yaw and roll angle in a tightening downward spiral \cite{stevens2015aircraft}.
The Dutch roll mode depends mainly on the \textit{dihedral} derivatives $C_{l_{\beta}}$ shown in \ref{eq:Cl} and it determines the amount of rolling in the Roll Subsidence mode. The coefficient matrices depend on the steady-state angle of attack and pitch attitude in both cases. Although they nominally apply to small perturbations about a steady-state flight condition, the equations can be used satisfactorily for perturbed roll angles of several degrees \cite{stevens2015aircraft}.   

\section{Optimal inputs for aircraft parameter estimation}\label{sec:OEDproblem}

In this section, an introduction to optimal maneuvers for aircraft parameter estimation is given and the formulation of an OED problem for the longitudinal and lateral aircraft dynamics is provided.

\subsection{Historical background and motivation}

In system identification in general as well as in the estimation of aerodynamic derivatives from flight tests, the design of the signal input provided for the system during experimental data collection is crucial for the accuracy of the subsequent parameter estimation. If a signal is not suitable for sufficient excitation of the relevant system dynamics, the data obtained during an experiment might not contain enough information on the desired parameters to allow for good estimation results. This creates need for a systematic design of optimal input signals for flight test maneuvers.

In the aerospace field, the importance of choosing appropriate control inputs for extraction of the aerodynamic derivatives from flight test data was first noticed by Gerlach \cite{gerlach1964analyse}. He proposed a qualitative method for the determination of optimal frequencies in scalar input signal to linear second order systems \cite{gerlach1971determination}. Important contribution to the theory and practice of the calculation of optimal aircraft input signals have been made subsequently by Mehra \cite{mehra1975status,mehra1974frequency,mehra1974time}. Based on the work of Kiefer, Wolfowitz \cite{kiefer1960equivalence} and Kiefer \cite{kiefer1959optimum}, Mehra proposed algorithms for the design of scalar and multi-dimensional input signals in the frequency domain as well as in time domain.    
An efficient method was implemented by Morelli where dynamics programming techniques were used to determine the optimal switching time of a input signal \cite{morelli1993oed}. The resulting input signals were of the \textit{bang-bang} type. Morelli's approach was afterwards applied by Cobleigh \cite{cobleigh1991design} and the resulting input signals were implemented by Noderer \cite{noderer1992analysis} for validation using real flight test data from an X-31 drop model. Nowadays optimal inputs are mainly designed in the aerospace field for \cite{licitra2017oed,klein1990optimal}:

\begin{itemize}
	\item reducing the number of expensive system identification flight tests,
	\item minimizing the length of flight test maneuver necessary to reach a specified level of accuracy of the aerodynamic derivatives,
	\item refinement and validation of the aerodynamic derivatives for control system analysis and design purposes and
	\item aircraft acceptance testing.
\end{itemize}

A fundamental problem in the design of input signals for parameter estimation is that the optimized design itself depends on the actual values of the unknown system parameters. As a consequence, these values would need to be known before the actual flight tests are made in order to optimize the experimental setup. However, if the parameters were already known, an estimation would obviously no longer be necessary. This problem is known as \textit{circularity problem} \cite{mulder1994identification}. Due to that, the optimal input design combined with parameter identification are used in practice in an iterative fashion, starting from a sufficient initial guess on the parameter values until a desired level of accuracy for the estimated values is met.

For the design of experiments for autonomous aircrafts which need to be performed in open-loop, additional considerations regarding trade-off between safety and accuracy of the aerodynamic derivatives need to be made.

\subsection{Optimum experimental design formulation}

Let us consider a mathematical model known a priori and defined as a set of ODE as in \ref{eq:ODE}

\begin{subequations}\label{eq:ODE}
	\begin{align}   
	\dot{\mathbf{x}}(t) & = \mathbf{f} \left(\mathbf{x}(t), \mathbf{u}(t), \mathbf{\theta}_{p}\right), \hspace{5pt} \mathbf{x}(0) = \mathbf{x}_{0} ,\hspace{5pt} t \in \left[0, T\right],  \\
	\mathbf{y}(t)       & = \mathbf{h} \left(\mathbf{x}(t), \mathbf{u}(t), \mathbf{\theta}_{p}\right), \\
	\mathbf{y}_{m}(i)   & = \mathbf{y}(i)   + \mathbf{\epsilon}(i),  \hspace{44pt}  i = 1,...,N 
	\end{align}
\end{subequations}
with differential states $\mathbf{x} \in \mathbb{R}^{n_{x}}$, output states $\mathbf{y} \in \mathbb{R}^{n_{y}}$, control inputs $\mathbf{u} \in \mathbb{R}^{n_{u}}$, a priori parameters $\mathbf{\theta}_{p} \in \mathbb{R}^{n_{\theta}}$. The output $\mathbf{y}_{m}$ is sampled in $N$ measurements along a time horizon $T$ and it is polluted by additive, zero-mean Gaussian noise \newline $\mathbf{\epsilon} \approx \eta(0, \mathbf{\Sigma_{y}})$ with $\mathbf{\Sigma_{y}} \in \mathbb{R}^{n_y \times n_y}$ the measurements noise covariance matrix. 

The Fisher information matrix $\mathbf{F}$ can then be expressed as 
\begin{equation}
\mathbf{F} = \sum_{i = 1}^{N} \left[ \left( \frac{\partial \mathbf{y}(i)}{\partial \mathbf{\theta}_{p}} \right)^{T} \mathbf{\Sigma_{y}}^{-1} \left( \frac{\partial \mathbf{y}(i)}{\partial \mathbf{\theta}_{p}} \right)\right].
\label{eq:FisherMatrix}
\end{equation}

The inverse of the Fisher information matrix $\mathbf{F}^{-1}$, which corresponds to the covariance matrix of the estimated parameters $\mathbf{\Sigma_{\theta}} \in \mathbb{R}^{n_{\theta} \times n_{\theta}}$, yields an universal lower bound on parameter estimation accuracies known as the \ac{CRLB}, which allows to optimize the input signal regardless the type of estimation algorithm implemented \cite{morelli1993oed}. The main idea of OED is to use an information function $\mathbf{\Psi}(\cdot)$ of $\mathbf{\Sigma}_{\theta}$ as the objective of an optimization problem. Therefore, a general model-based OED problem which considers input and output constraints can be formulated as
\begin{subequations}\label{eq:OCP_TC}
	\begin{align}
	\underset{\mathbf{x}(\cdot),\mathbf{u}(\cdot)}{\text{minimize}} \hspace{1cm} & \mathbf{\Psi} \left(\mathbf{\Sigma_{\theta}} \left[\mathbf{x}(\cdot), \mathbf{u}(\cdot), \mathbf{\theta}_{p}\right] \right) \label{eq:OCP_objective} \\ 
	\text{subject to:} \hspace{1cm} & \dot{\mathbf{x}}(t) = \mathbf{f} \left(\mathbf{x}(t), \mathbf{u}(t), \mathbf{\theta}_{p}\right), \hspace{0.14cm}  t \in \left [ 0,T  \right ],\label{eq:OCP_ode}  \\ 
	& \mathbf{x}(0) = \mathbf{x}_{0}, \\
	& \mathbf{u}_{\mathrm{min}} \leq \mathbf{u}(t) \leq \mathbf{u}_{\mathrm{max}},  \hspace{0.7cm}  t \in  \left [ 0,T  \right ], \label{eq:OCP_control_bounds} \\
	& \mathbf{x}_{\mathrm{min}} \leq \mathbf{x}(t)  \leq \mathbf{x}_{\mathrm{max}}, \hspace{0.75cm}  t \in \left [ 0,T  \right ]. \label{eq:OCP_state_bounds}
	\end{align}
\end{subequations}
Different information functions can be used in the optimization problem \ref{eq:OCP_TC} with different features \cite{gupta1975input,nahi1969optimal,nahi1971design}.We will further use the so called \textit{A-criterion} $\mathbf{\Psi_{A}}(\cdot)$ \cite{mehra1975status} within the objective of the OED problem, which is the scaled trace of $\mathbf{\Sigma_{\theta}}$ as in \ref{eq:Acriterion}
\begin{equation}\label{eq:Acriterion}
\mathbf{\Psi_{A}} \left( \mathbf{\Sigma_{\theta}} \right) = \frac{1}{n_{\theta}} \cdot \mathrm{trace} \left( \mathbf{\Sigma_{\theta}} \right) =  \frac{1}{n_{\theta}} \sum_{i=1}^{n_\theta} \mathbf{\Sigma}_{\mathbf{\theta},(i,i)} = \frac{1}{n_{\theta}} \sum_{i=1}^{n_\theta} \mathrm{Var}(\mathbf{\theta}_{p,(i)}),
\end{equation}
so that using $\mathbf{\Psi_{A}}(\cdot)$, one optimize the experimental setup in terms of minimizing the sum of the variances of the unknown parameters.

\subsection{Algorithm implementation}

Within this work, the optimum experimental designs are computed using \textsc{casiopeia}, an open-source tool for parameter estimation and OED \cite{Buerger2016} based on \textsc{CasADi} \cite{Andersson2013b}. \textsc{casiopeia} computes the covariance matrix $\mathbf{\Sigma_{\theta}}$ from the inverse of the KKT Matrix of the underlying parameter estimation problem using a Schur complement approach. Details on method and implementation can be found in \cite{Buerger2017}.

Provided the system dynamics \ref{eq:OCP_ode}, a discretization time grid, bound specifications for variables as in \ref{eq:OCP_control_bounds} \ref{eq:OCP_state_bounds} and an initial guess for the parameter values $\mathbf{\theta}_\mathrm{init}$ and for the input signal $\mathbf{u}_\mathrm{init}$, the continuous-time optimization problem is discretized and formulated as a NLP automatically by \textsc{casiopeia} using direct collocation \cite{Biegler2010} with Lagrange polynomials. The resulting NLP is solved using \textsc{IPOPT} \cite{Waechter2006} with linear solver \textsc{MA86} \cite{HSL2017} to obtain improved input signals $\mathbf{u}_\mathrm{opt}$.

If for the initial values $\mathbf{\theta}_\mathrm{init}$ of the parameters used within OED it holds for two parameters $\theta_i, \theta_j \in \mathbf{\theta}$ with $i\neq j$ that $\mathbf{\theta}_{\mathrm{init},i} > \mathbf{\theta}_{\mathrm{init},j}$, it is likely that $\mathrm{Var}(\theta_i) > \mathrm{Var}(\theta_j)$. Due to the higher contribution of $\mathrm{Var}(\theta_i)$ to \ref{eq:OCP_objective}, the optimizer then might overly increase certainty of $\theta_i$ while disregarding to increase or even decreasing certainty of $\theta_j$. Due to that, $\mathbf{\theta}$ is within the OED problem formulations of this work not introduced by the values of $\mathbf{\theta}_\mathrm{init}$, but as a vector of entries 1 scaled by the corresponding entries of $\mathbf{\theta}_\mathrm{init}$ to reduce the effects of the numerical values of $\mathbf{\theta}_\mathrm{init}$ on the OED result.

\subsection{A priori model}

The OED problem \ref{eq:OCP_TC} requires an a priori model with $\mathbf{\theta}_\mathrm{init}$ of sufficient accuracy. Various methods can be applied to obtain a priori models. If the airframe is similar to an existing aircraft, its model can be scaled. For instance, the Digital DATCOM \cite{hoak1975usaf} is a purely empirical guide to estimating aerodynamic derivatives based on aircraft configuration and the experience of engineers. If the airfoils and aircraft configuration are new, one can perform analysis via the lifting line method \cite{versteeg2007introduction}, \ac{CFD} \cite{anderson2017fundamentals}, wind-tunnel tests or previous flight tests. Depending on the available resources, combinations of these methods can be used.

In this work, a priori models are retrieved from both the lifting line method and previous flight tests \cite{licitra2017pe}. A steady wing-level flight condition was obtained for $\VTe = 20\,\mathrm{m/s}$ and flap setting $\delta_{f} = 0 \,^\circ$. The equilibrium point is hold for $\delta_{e} = -1.5 \, ^\circ, \alpha_{\mathrm{e}} = - 0.4 \, ^\circ$ and $\theta_{\mathrm{e}} = - 4.5 \, ^\circ$ with the other states equal to zero.

The \ac{OED} problem \ref{eq:OCP_TC} takes into account longitudinal \ref{eq:longDyn} and lateral dynamics \ref{eq:latDyn} separately with a priori dimensional aerodynamic derivatives shown in tables~\ref{tab:LongitudinalDerivatives},\ref{tab:LateralDerivatives}. Note that $M_{V}$ and $Y_{p}$ were not included in $\theta_\mathrm{init}$ since they were zero. Yet, the LTI systems \ref{eq:longDyn},\ref{eq:latDyn} were augmented with their input derivatives in order to take into account the rate of control surfaces. 

As mentioned in section~\ref{sec:steadyConditionAndDecoupling}, these a priori models also provide an insight into the general characteristics of the aircraft behavior via modal analysis. In table~\ref{tab:AircraftMode} aircraft modes with the corresponding natural frequencies $\omega_{n}$, damping ratios $\delta$, constant times $\tau$, overshoots $S_{\%}$ and period of oscillations $P_{O}$ are shown.

\begin{table}[tbhp]
	\caption{Aircraft modal analysis for steady wing-level flight at $V_{T_{e}} = 20 \mathrm{[m/s]}$ }
	\label{tab:AircraftMode}
	\centering
	\begin{tabular}{lccccc}\hline
		\bf Mode        &$\omega_{n}\;\mathrm{[rad/s]}$&$\delta\;\mathrm{[-]}$&$\tau\;\mathrm{[s]}$&$S_{\%}\;\mathrm{[\%]}$&$P_{O}\;\mathrm{[s]}$\\\hline
		Phugoid         &         0.52                 &       0.09           &      1.94          &       74.06       &      12.23       \\
		Short-period    &         3.72                 &       0.84           &      0.27          &        0.83       &       3.08       \\
		Spiral          &         0.08                 &       - -            &      11.74         &        - -        &       - -        \\
		Dutch roll      &         2.09                 &       0.21           &      0.48          &        50.55      &       3.08       \\
		Roll Subsidence &         11.12                &       1.0            &      0.09          &        - -        &       - -        \\\hline	
	\end{tabular}
\end{table} 

The modal analysis shows that both Phugoid and Spiral mode are located below $0.1\,\mathrm{Hz}$, hence the aerodynamic derivatives relative to these modes will not be identified accurately due to a lack of excitation in this frequency range. However, this is not of great concern, as the slow nature of these dynamics can be easily handled by a pilot or a control system \cite{dorobantu2013system}. Finally, the Spiral mode is lightly unstable as expected.

\subsection{Constraints selection}
In practice, it is hardly possible to apply input signals which correspond to full deflection of the control surfaces without exceeding the limits of the permissible flight envelope. On the one hand, one must scaled down the input signal amplitude in order to restrict the aircraft response within a region for which the model structure assumed in \ref{eq:OCP_TC} is valid. Note that, besides the linear structures used in this work, the aerodynamic properties described in section~\ref{sec:aerodynamics} are only valid for a neighborhood around a specific steady condition. On the other hand, if input signals are scaled up, then the estimation accuracy is enhanced due to a higher \ac{SNR}. 

Therefore, constraints in \ref{eq:OCP_control_bounds},\ref{eq:OCP_state_bounds} should be enforced in order to ensure the system response close to a specific steady wing-level flight condition without any flight envelope violation and at the same time guarantee an acceptable \ac{SNR}. In this work, constraints for the OED problem \ref{eq:OCP_TC} have been chosen as follows:
\begin{itemize}
	\item control surface deflections $(\delta_{\mathrm{a}},\delta_{\mathrm{e}},\delta_{\mathrm{e}})$, angle of attack $\alpha$ as well as the airspeed $\VT$ were constrained in order to keep the aircraft within the region where the linear model is applicable;
	\item rate of control surface deflections $(\dot{\delta_{\mathrm{a}}},\dot{\delta_{\mathrm{e}}},\dot{\delta_{\mathrm{r}}})$ were constrained according to the maximum speed of the installed servos;
	\item the body angular rates $(p,q,r)$ and Euler angles $(\phi,\theta,\psi)$ should be bounded with respect to the flight envelope limits since any violation of the flight envelope would result in abortion of the system identification test. To account for model mismatch and inaccuracies of the a priori model, these bounds were enforced with a safety margin of $20\,\%$ w.r.t. the flight envelope limits.    
\end{itemize} 
Table~\ref{tab:FlightEnvelope} summarizes the flight envelope, input and state constraints relative. 
\begin{table}[tbhp]
	\caption{Flight envelope and OED constraints relative to the trim condition at $\VTe = 20 \mathrm{[m/s]}$}
	\label{tab:FlightEnvelope}
	\centering
	\begin{tabular}{ccccc}\hline
		\bf Variable & \bf Flight envelope   &\bf OED constraints &\bf Range & \bf Unit           \\\hline
		$V_{T}$      &     ( 12,30)          &       ( 17,23)     &     6    & $[\mathrm{m/s}]$   \\
		$\beta$      &     (-20,20)          &      (-7.5,7.5)    &    15    & $[\mathrm{deg}]$   \\
		$\alpha$     &     (-8,20)           &     (-4.36,3.64)   &     8    & $[\mathrm{deg}]$   \\ 
		$\phi$	     &     (-35,35)          &       (-28,28)     &    56    & $[\mathrm{deg}]$   \\
		$\theta$     &     (-30,40)          &    (-28.77,27.33)  &    56    & $[\mathrm{deg}]$   \\
		$p$          &     (-60,60)          &       (-48,48)     &    96    & $[\mathrm{deg/s}]$ \\
		$q$          &     (-40,40)          &       (-32,32)     &    64    & $[\mathrm{deg/s}]$ \\
		$r$          &     (-40,40)          &       (-32,32)     &    64    & $[\mathrm{deg/s}]$ \\
		$(\delta_{a},\delta_{e},\delta_{r})$ &          (-,-)     &     (-5,5)      & 10    & $[\mathrm{deg}]$\\
		$(\dot{\delta_{a}},\dot{\delta_{e}},\dot{\delta_{r}})$    &       (-,-)     &     (-3.25,3.25)      & 7    & $[\mathrm{rad/s}]$\\\hline		
	\end{tabular}
\end{table}

\subsection{Control input initialization}
\begin{figure}[tbhp]
	\centering
	\includegraphics[scale = 0.48]{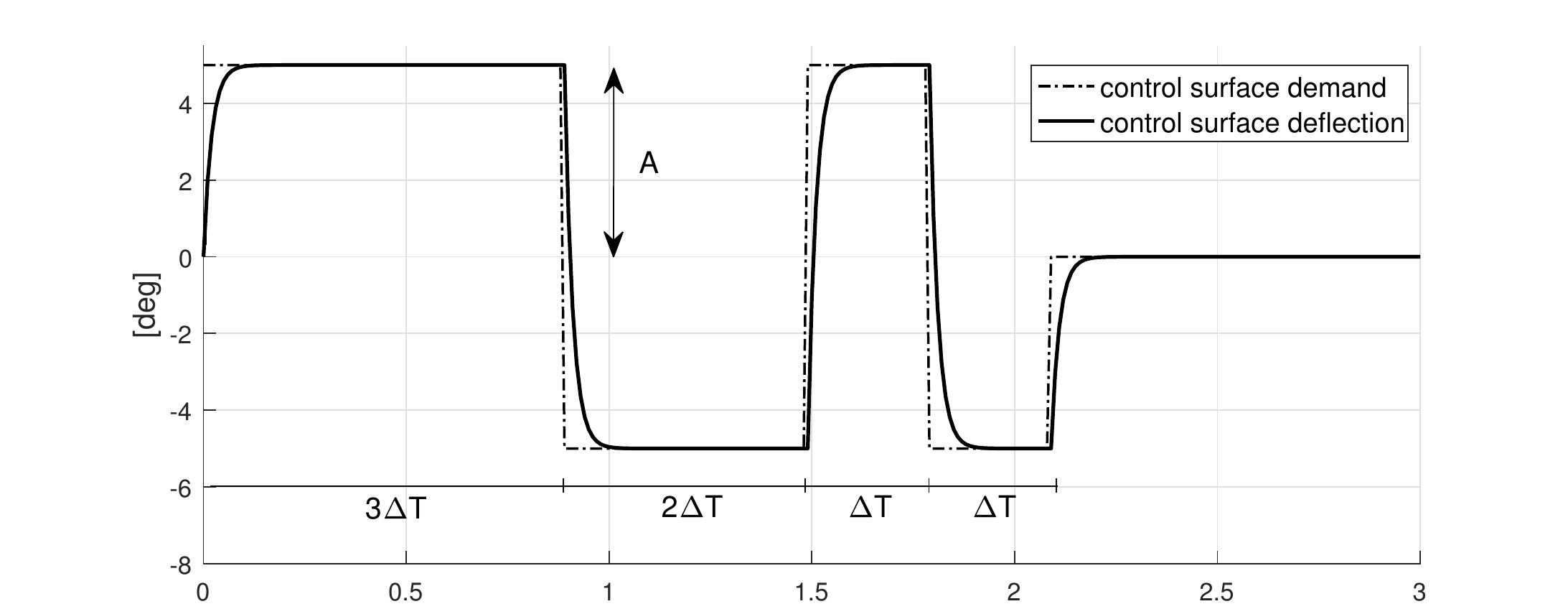}
	\caption{\textit{3-2-1-1 maneuver} of amplitude $A$ and duration $7\Delta T$. In dash-dot line the control surface demand delivered by the FCC and in solid line the actual control surface deflection.}
	\label{fig:3-2-1-1maneuver}
\end{figure}

The optimization problem \ref{eq:OCP_TC} needs to be initialized with a suitable, initial input signal $\mathbf{u}_\mathrm{init}$. The initial control input chosen for this application is widely used in the aerospace field due to its easy implementation and good estimation performance. Such maneuver comes from an optimization procedure of a sequence of step functions, developed by Koehler \cite{koehler1977auslegung}. The aim of Koehler was to find a signal with a shape as simple as possible and power distributed uniformly over a wide range of frequencies \cite{plaetschke1979practical}. The input signal has a bang-bang behavior with a duration $7\Delta T$ with switching times at $t = 3\Delta T$, $t = 5\Delta T$, and $t = 6\Delta T$. For this reason, such an input signal is called a \textit{3-2-1-1 maneuver}, see fig.~\ref{fig:3-2-1-1maneuver}.

In \cite{mulder1994identification}, it was shown that the 3-2-1-1 maneuver provides the best estimation accuracy for both aircraft longitudinal and lateral dynamics among \textit{Doublets}, \textit{Mehra}, \textit{Schulz} and \textit{DUT} input signals. Yet, it is shown that only Doublets and 3-2-1-1 input signals contain significant amount of power above $1\,\mathrm{Hz}$, though the 3-2-1-1 maneuver embraces much higher frequencies compared to Doublets.     

\section{Results of the experimental design}\label{sec:Results}
In this section, the results of the experimental design are shown and assessed with respect to the conventional 3-2-1-1 maneuvers.

\subsection{Optimization of the experimental design}\label{sec:SimulationResult}
Three OED problems \ref{eq:OCP_TC} are computed, one for the longitudinal dynamics \ref{eq:longDyn} and two for the lateral dynamics \ref{eq:latDyn}. Since only one single axis can be excited at a time,  when the aileron is chosen as input, the rudder is assumed zero along the entire experiment. As a consequence, the aerodynamic derivatives relative to the rudder ($Y_{\delta_{\mathrm{r}}},L_{\delta_{\mathrm{r}}}',N_{\delta_{\mathrm{r}}}'$) will be \textit{structurally unidentifiable} with $\mathbf{F}$ not of full rank. In this case, one has to discard the corresponding unidentifiable parameters to prevent rank deficiency of $\mathbf{F}$ \cite{mulder1994identification}. Same considerations are valid in the case when the rudder is used as input and the aileron is kept zero. 
Note that in practice, when the roll axis is excited by aileron deflection, the rudder stabilizes the yaw axis (and the elevator the pitch axis) during the entire system identification experiment. Hence, the rudder control surface will slightly differ from zero. For the sake of comparability, the amplitude of $\mathbf{u}_\mathrm{init}$ and $\mathbf{u}_\mathrm{opt}$ are set to equal value. 

The 3-2-1-1 maneuvers are chosen through both a qualitative consideration in the frequency domain \cite{marchand1977untersuchung} and a \textit{trial-and-error} approach in order to ensure that the system response is within the prescribed constraints.

\begin{figure}[tbhp]
	\centering
	\includegraphics[width = 370pt, height = 300pt]{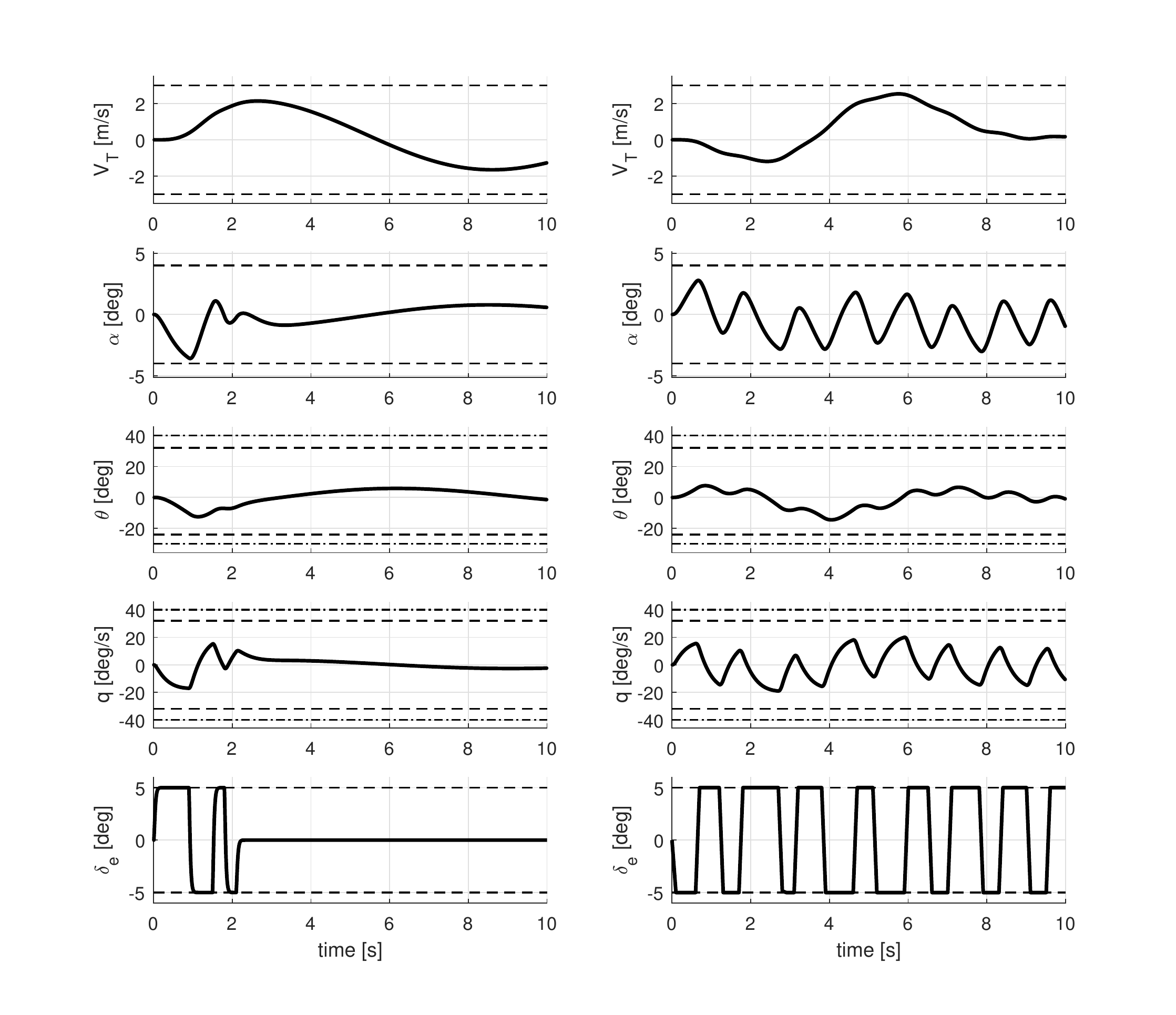}
	\caption{Longitudinal response using the 3-2-1-1 maneuver (left column) and optimal response (right column). Dash line the OED constraints whereas dash-dot line the flight envelope.}
	\label{fig:LongDyn_initVSopti}
\end{figure}

\begin{figure}[tbhp]
	\centering
	\includegraphics[width = 370pt, height = 300pt]{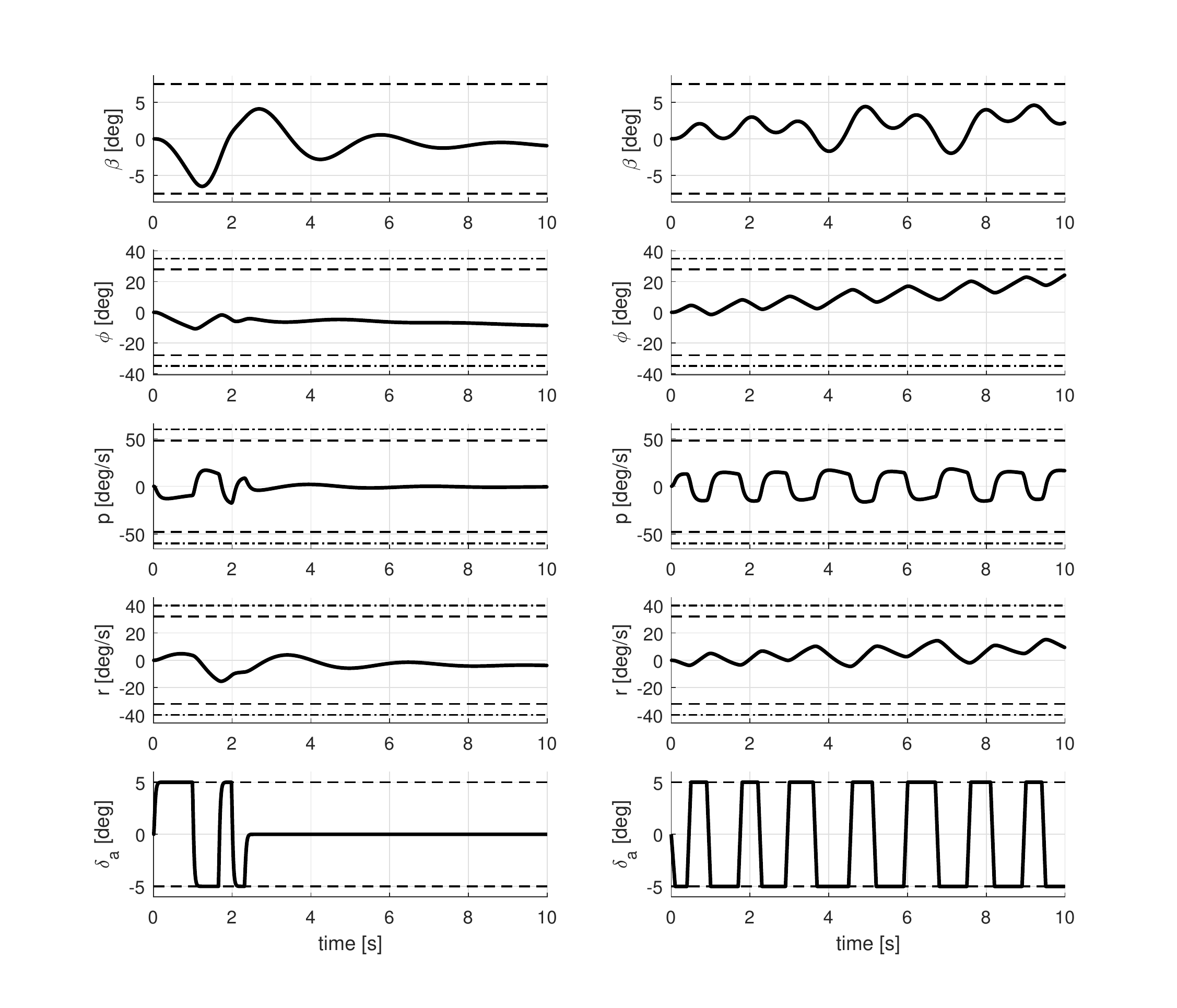}
	\caption{Lateral response via aileron deflection using the 3-2-1-1 maneuver (left column) and optimal response (right column). Dash line the OED constraints whereas dash-dot line the flight envelope.}
	\label{fig:LatDyn_initVSopti_da}
\end{figure}

\begin{figure}[tbhp]
	\centering
	\includegraphics[width = 370pt, height = 300pt]{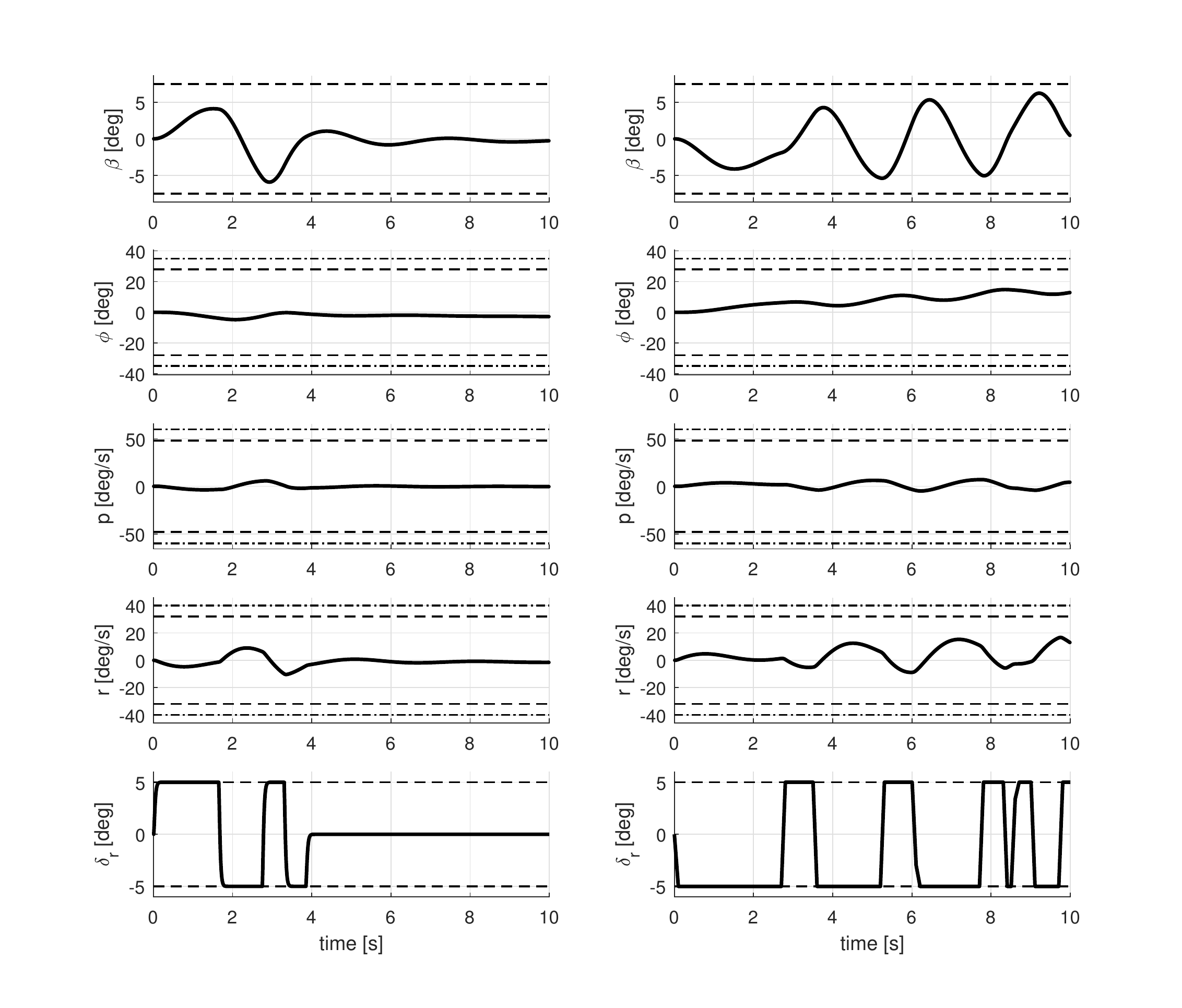} 
	\caption{Lateral response via rudder deflection using the 3-2-1-1 maneuver (left column) and optimal response (right column). Dash line the OED constraints whereas dash-dot line the flight envelope.}
	\label{fig:LatDyn_initVSopti_dr}
\end{figure}

The responses of the a priori LTI systems for the optimized inputs obtained from the solution of \ref{eq:OCP_TC} are shown in figs.~ \ref{fig:LongDyn_initVSopti},\ref{fig:LatDyn_initVSopti_da},\ref{fig:LatDyn_initVSopti_dr} with the corresponding responses to the 3-2-1-1 maneuvers.

It turns out that the OED problem \ref{eq:OCP_TC} leads to a signal of \textit{bang-bang} type. Such outcome is in agreement with analytic results \cite{chen1975input} and previous flight test evaluations which demonstrate that square wave type inputs are superior to sinusoidal type inputs for parameter estimation experiments, largely due to their wider frequency spectrum \cite{mulder1986design}. More precisely, the signal inputs resemble modulated square waves with a finite slope of $2.5 \, \mathrm{rad/s}$ due to the rate of deflection constraints. 

In fig.~\ref{fig:LongDyn_initVSopti} one can observe the decoupling between the Phugoid mode which dominates the airspeed $\VT$ and pitch $\theta$ responses, with the fast changes on the angle of attack $\alpha$ and pitch rate $q$ coming from the Short-period mode. 

As shown in fig.~\ref{fig:LatDyn_initVSopti_da}, the optimal lateral response caused by the aileron deflection shows a good excitation on the roll rate $p$. The cross-coupling involves a modest excitation on the yaw rate $r$ and side-slip angle $\beta$ whereas the roll angle $\phi$ drifts slowly towards the edge of the admissible range due to the unstable \textit{Spiral} mode.

Regarding the yaw excitation via rudder deflection depicted in fig.~\ref{fig:LatDyn_initVSopti_dr}, an optimal response would be provided by setting the rudder at the maximum allowable deflection for approximately $3\,\mathrm{s}$ such that an oscillatory motion with a gradual increment on amplitude on the side slip angle $\beta$ as well as the roll rate $r$ is triggered. Subsequently, a bang-bang behavior is carried out so as to avoid constraint violations. The roll rate $p$ is barely excited due to cross-coupling though, roll angle $\phi$ slowly diverges as in the previous case.               


\subsection{Performance assessments}\label{sec:PerformanceAssessments}
In this section, optimal maneuvers are assessed by the \ac{CRLB} which is the theoretical lower limits for parameter standard errors using an efficient and asymptotically unbiased estimator, such as maximum likelihood. The \ac{CRLB} depends on the diagonal entries of the Fisher information matrix $\mathbf{F}$ \ref{eq:FisherMatrix} which is formally \cite{soderstrom2001identification,tischler2006aircraft}
\begin{equation}\label{eq:std-CR}
\sigma_{i} \geq \mathrm{CRLB}_{i} = \frac{1}{\sqrt{\mathbf{F}_{ii}}}.
\end{equation}
A performance analysis of signal inputs computed via the \ac{CRLB} isolates the merits of the input design from the merits of the parameter estimation algorithm used to extract the aerodynamic derivatives from the flight data \cite{klein1990optimal}. Yet, the relation between the parametric uncertainty and \ac{CRLB} allows to form a comprehensive uncertain aircraft model. 

Several factor can cause high \ac{CRLB} values, e.\,g.: 
\begin{itemize}
	\item from a optimization point of view, large values for $\mathrm{CRLB}_{i}$ denote a low curvature in the cost function, i.\,e., a high insensitivity with respect to the $i^\mathrm{th}$ parameter \cite{tischler2006aircraft};
	\item from a system identification point of view, high \ac{CRLB} values indicate either that the $i^\mathrm{th}$ parameter is physically insignificant with respect to the measured aircraft response or that there exists a correlation between parameters, i.\,e., these parameters can vary together, making their individual values difficult to determine \cite{soderstrom2001identification}. 
\end{itemize}
Tables~\ref{tab:CRlongDyn},\ref{tab:CRlatDyn_da},\ref{tab:CRlatDyn_dr} show the \ac{CRLB} values for the optimal system responses ($\mathrm{CRLB_{opt}}$) with the corresponding initial responses ($\mathrm{CRLB_{init}}$) whereas $\mathrm{\Delta CRLB_{\%}}$ indicates the percent deviation between $\mathrm{CRLB_{opt}}$ and $\mathrm{CRLB_{init}}$. 
A negative value for $\mathrm{\Delta CRLB_{\%}}$ indicates an improvement in terms of estimation accuracy for the $i^\mathrm{th}$ parameter, and vice versa for positive values.   
\begin{table}[tbhp]
	\caption{dimensional aerodynamic derivatives longitudinal dynamics}
	\label{tab:CRlongDyn}
	\centering
	\begin{tabular}{crrrc}\hline
		\bf Derivatives  &\bf Value &$\mathrm{CRLB_{init}}$ &$\mathrm{CRLB_{opt}}$ &$\mathrm{\Delta CRLB_{\%}}$\\\hline
		$X_{V}$               &  -0.147  &   0.1978    &  0.1147   &  -42.0  \\
		$X_{\alpha}$          &   7.920  &  14.0706    &  8.6891   &  -38.2  \\
		$X_{q}$               &  -0.163  &   2.9936    &  1.5916   &  -46.8  \\			
		$X_{\delta_{e}}$      &  -0.232  &   3.4318    &  1.9859   &  -42.1  \\			
		$Z_{V}$               &  -0.060  &   0.0007    &  0.0004   &  -42.0  \\
		$Z_{\alpha}/\VTe$     &  -4.400  &   0.0491    &  0.0303   &  -38.2  \\
		$Z_{q}$               &   0.896  &   0.0104    &  0.0056   &  -46.8  \\			                    
		$Z_{\delta_{e}}/\VTe$ &  -0.283  &   0.0120    &  0.0069   &  -42.1   \\
		$M_{\alpha}$          &  -6.180  &   0.0098    &  0.0061   &  -38.2   \\
		$M_{q}$               &  -1.767  &   0.0021    &  0.0011   &  -46.8   \\	   
		$M_{\delta_{e}}$      & -10.668  &   0.0024    &  0.0014   &  -42.1   \\\hline		
	\end{tabular}
\end{table}
\begin{table}[tbhp]
	\caption{dimensional aerodynamic derivatives lateral dynamics: aileron input}
	\label{tab:CRlatDyn_da}
	\centering
	\begin{tabular}{crrrc}\hline
		\bf Derivatives  &\bf Value &$\mathrm{CRLB_{init}}$ &$\mathrm{CRLB_{opt}}$ &$\mathrm{\Delta CRLB_{\%}}$\\\hline
		$Y_{\beta}/\VTe$      &  -0.167  &   0.8907    &  1.1087    &   24.5       \\
		$ Y_{r}   $           &  -0.976  &   0.3825    &  0.3526    &   -7.8       \\
		$Y_{\delta_{a}}/\VTe$ &  -0.046  &   1.7305    &  0.5869    &  -66.1       \\	  
		$L_{\beta}'$          &  -8.201  &   0.8907    &  1.1087    &   24.5       \\
		$L_{p}'$              & -11.292  &   0.6856    &  0.1504    &  -78.1       \\
		$L_{r}'$              &   3.853  &   0.3825    &  0.3526    &   -7.8       \\
		$L_{\delta_{a}}'$     & -32.600  &   1.7305    &  0.5869    &  -66.1       \\				
		$N_{\beta}'$          &   3.214  &   0.8907    &  1.1087    &   24.5       \\			
		$N_{p}'$              &  -0.750  &   0.6856    &  0.1504    &  -78.1       \\
		$N_{r}'$              &  -0.457  &   0.3825    &  0.3526    &   -7.8       \\			    
		$N_{\delta_{a}}'$     &   0.716  &   1.7305    &  0.5869    &  -66.1       \\\hline		
	\end{tabular}
\end{table}
\begin{table}[tbhp]
	\caption{dimensional aerodynamic derivatives lateral dynamics: rudder input}
	\label{tab:CRlatDyn_dr}
	\centering
	\begin{tabular}{crrrc}\hline
		\bf Derivatives  &\bf Value &$\mathrm{CRLB_{init}}$ &$\mathrm{CRLB_{opt}}$ &$\mathrm{\Delta CRLB_{\%}}$\\\hline
		$Y_{\beta}/\VTe$      & -0.167  &  10.0229    & 3.0396     & -69.7    \\
		$ Y_{r}   $           & -0.976  &   3.3870    & 0.9643     & -71.5    \\
		$Y_{\delta_{r}}/\VTe$ &  0.093  &   1.3518    & 0.5674     & -58.0    \\	  
		$L_{\beta}'$          & -8.201  &  10.0229    & 3.0396     & -69.6    \\
		$L_{p}'$              &-11.292  &  12.2297    & 3.7003     & -69.7    \\			
		$L_{r}'$              &  3.853  &   3.3870    & 0.9643     & -71.5    \\
		$L_{\delta_{r}}'$     &  0.524  &   1.3518    & 0.5674     & -58.0	  \\
		$N_{\beta}'$          &  3.214  &  10.0229    & 3.0396     & -69.6    \\					
		$N_{p}'$              & -0.750  &  12.2297    & 3.7003     & -69.7    \\
		$N_{r}'$              & -0.457  &   3.3870    & 0.9643     & -71.5    \\			    
		$N_{\delta_{r}}'$     & -2.370  &   1.3518    & 0.5674     & -58.0    \\\hline		
	\end{tabular}
\end{table}

The optimal signal input for the longitudinal dynamics (table~ \ref{tab:CRlongDyn}, fig.~\ref{fig:LongDyn_initVSopti}) provides an overall increment of $\approx 40\%$ in terms of estimation accuracy. Though, as expected, the parameters determined by the Phugoid mode, i.\,e., $X_{V},X_{\alpha},X_{q}$ and $X_{\delta_{e}}$, are still subject to high uncertainty due to the time window of the experiment set to $10\, \mathrm{s}$. Yet, the 3-2-1-1 maneuver used as initial signal input provides more than the acceptable accuracy for  $Z_{V},Z_{\alpha},Z_{q},Z_{\delta_{e}},M_{\alpha},M_{q}$ and $M_{\delta_{e}}$ which are derivatives related to the Short-period mode.  

The optimal lateral response via aileron deflection (table~\ref{tab:CRlatDyn_da},fig.~\ref{fig:LatDyn_initVSopti_da}), reduce mainly the uncertainty on the approximation of the aileron to roll rate transfer function which is  
$\frac{p\left(s\right)}{\delta_{a}\left(s\right)} = \frac{L_{\delta_{a}}'}{s - L_{p}'}$ \cite{stevens2015aircraft}.
The contribution of $Y_{\delta_{a}}$ with respect to the overall aircraft response appears negligible for this steady configuration, hence its uncertainty will be high in the optimal case, too.  
On the other hand, significant increment in terms of accuracy is shown in $N_{p}'$ and $N_{\delta_{a}}'$ which are parameters relative to the yaw moment due to the roll rate $p$ and aileron deflection $\delta_{\mathrm{a}}$, respectively. One can observe that remain derivatives are either poorly improved ($Y_{r},L_{r}',N_{r}'$) or they experience a loss of accuracy ($Y_{\beta},L_{\beta}',N_{\beta}'$). This is not surprising since these derivatives are all related to the yaw dynamics which is barely excited during aileron deflection.

Finally, table~\ref{tab:CRlatDyn_dr} shows that the optimal rudder deflection in fig.~\ref{fig:LatDyn_initVSopti_dr} provides a meaningful improvement mainly on parameters relative to the yaw dynamics.  

\section{Application of the optimal inputs within real flight tests}\label{sec:RealImplementation}
Within this section, the optimized maneuvers are first validated via reliable flight simulator for safety issues and subsequently experimental data coming from a real flight test campaign are shown.

\subsection{Signal inputs set-up and safety assessments}
As already mentioned, the case study is an autonomous system hence, a flight plan must be set. The \ac{FCC} of the case study allows to define control surfaces demands as steps with tunable amplitude and time length only. Therefore, the steps transition shown in Figs.\ref{fig:LongDyn_initVSopti},\ref{fig:LatDyn_initVSopti_da},\ref{fig:LatDyn_initVSopti_dr} are approximated as tight step functions as in fig.~\ref{fig:3-2-1-1maneuver}. 
Before that any real system identification flight test can be performed, each signal input need to be validated via reliable flight simulator for different wind conditions as well as degree of parameter's uncertainty so as to get a better confidence relative to the safety of the flight test.

\begin{figure}[tbhp]
	\centering
	\includegraphics[scale = 0.50]{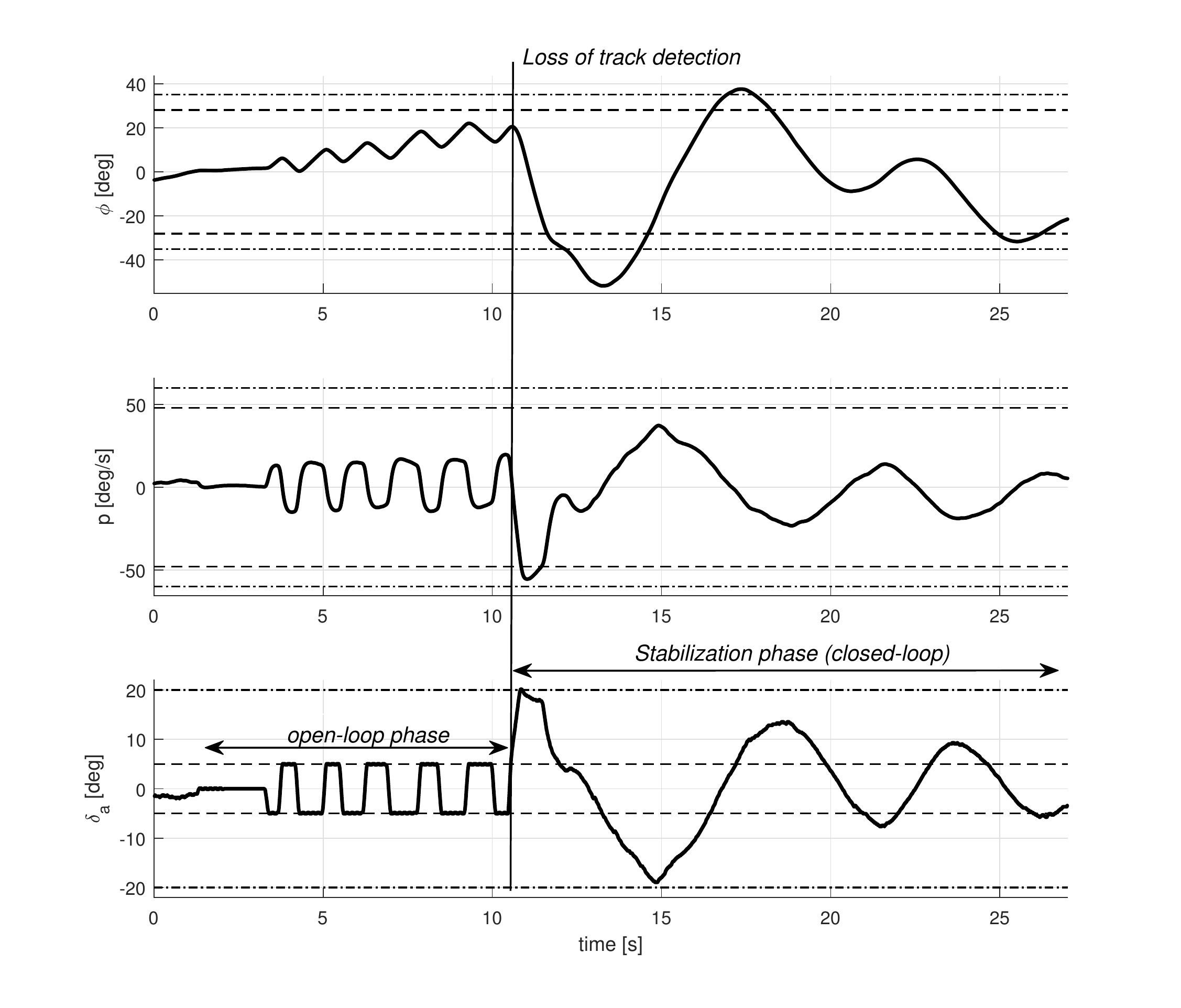}
	\caption{Simulated flight test during roll excitation via optimal aileron deflection. The sequence start setting the aircraft at steady state wing level trim condition; afterwards, the optimal signal input is performed along the roll axis via $\delta_{a}$. The roll angle $\phi$ oscillates and move far form the trim condition which causes a deviation from the flight path. At $\approx 11 \mathrm{[s]}$ flight envelope protection triggers and the open-loop sequence stops, recovering the aircraft attitude via feedback controls. In dash lines the OED constraints whereas dash-dot line the flight envelope for $\phi,p$ and limiters from $\delta_{a}$.}
	\label{fig:AbortOptimalAileronGraph}
\end{figure}

\begin{figure}[tbhp]
	\centering
	\includegraphics[scale = 0.35]{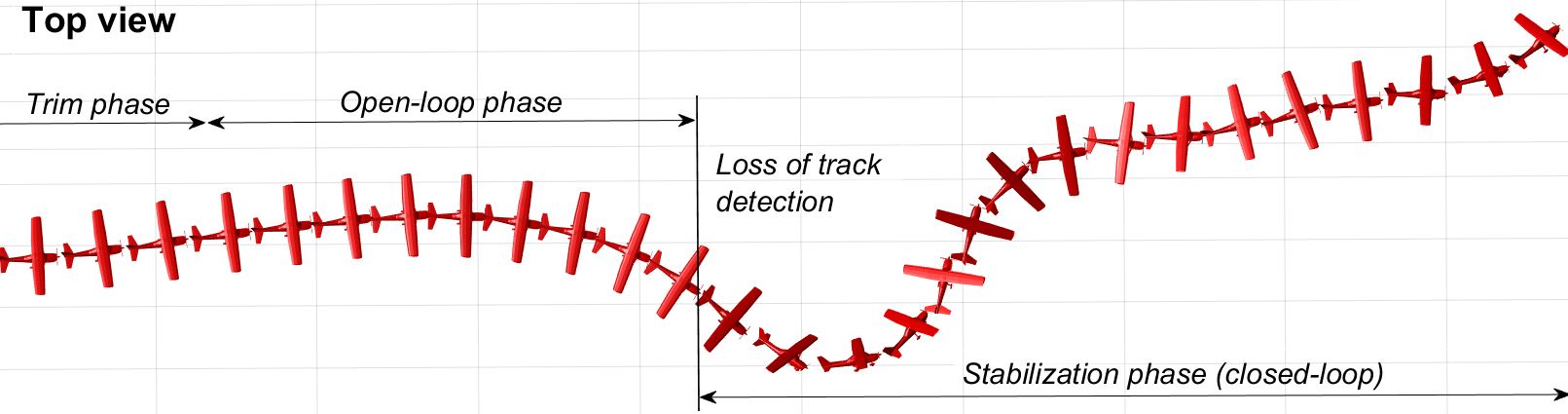}
	\caption{Top view of the aircraft during loss of track caused by the optimized aileron sequence}
	\label{fig:AbortOptimalAileron3D}
\end{figure}

The results obtained from the high fidelity simulator designed by Ampyx Power B.V. \cite{AP} has shown that the aircraft was able to complete successfully the system identification flight test with optimal elevator deflection and rudder without any flight envelope violation, providing good excitation of the longitudinal and yaw dynamics, respectively.  
As far as it regards the excitation of the lateral dynamics via aileron deflection, it turns out that the vehicle is prone to loss of the flight path due to cross-coupling effect and unstable Spiral mode (see figs.\ref{fig:AbortOptimalAileronGraph},\ref{fig:AbortOptimalAileron3D}).   
Therefore, for safety reasons part of the optimal aileron sequence is discarded and propulsion system are not turned off in order to avoid significant deviation from the steady condition. Note that, the thrust effect provides negligible model mismatch on the roll axis. Alternatively, one could superimpose the pilot command with the actual optimal signal input so as to prevent loss of track \cite{dorobantu2013system}, though degradation of the estimation accuracy might occurs. 

\subsection{Experimental results}
A real flight test campaign is conducted using the optimized inputs obtained in section \ref{sec:SimulationResult}.
During the experiments, the estimated with speed was $\approx 2 \, \mathrm{m/s}$ hence, low process noise caused by turbulences is expected. Each signal input is performed multiple times so has to obtain a trend of the aircraft responses.

\begin{figure}[tbhp]
	\centering
	\includegraphics[width = 370pt, height = 300pt]{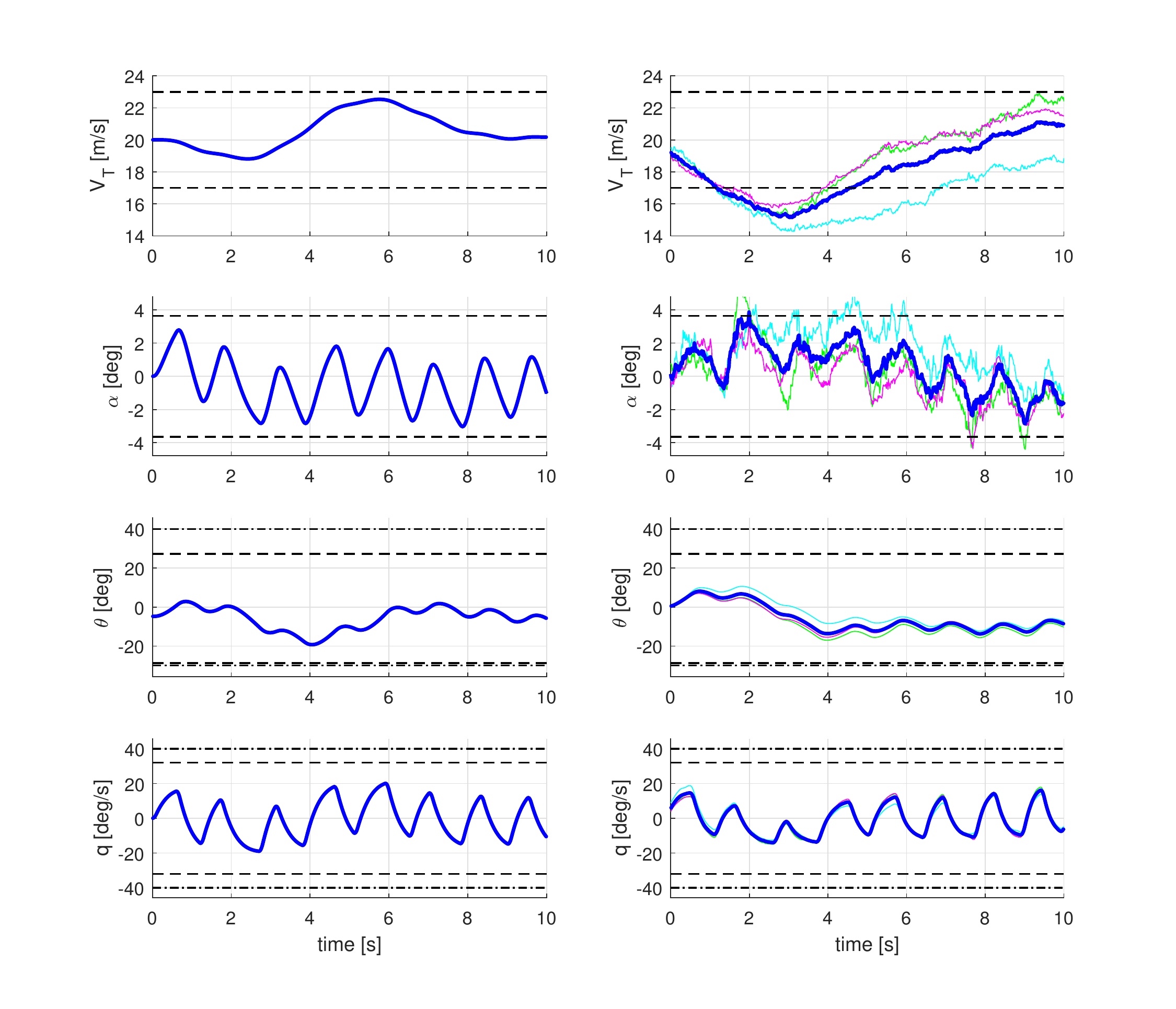}
	\caption{Simulated optimal longitudinal response of the a priori model (left column) versus three real experiments (right column). The real experiments are in green, magenta and cyan thin solid lines whereas their average is in thick blue solid line. In dash line the OED constraints while dash-dot line the flight envelope.}
	\label{fig:LongDyn_simVSreal}
\end{figure}

\begin{figure}[tbhp]
	\centering
	\includegraphics[width = 370pt, height = 300pt]{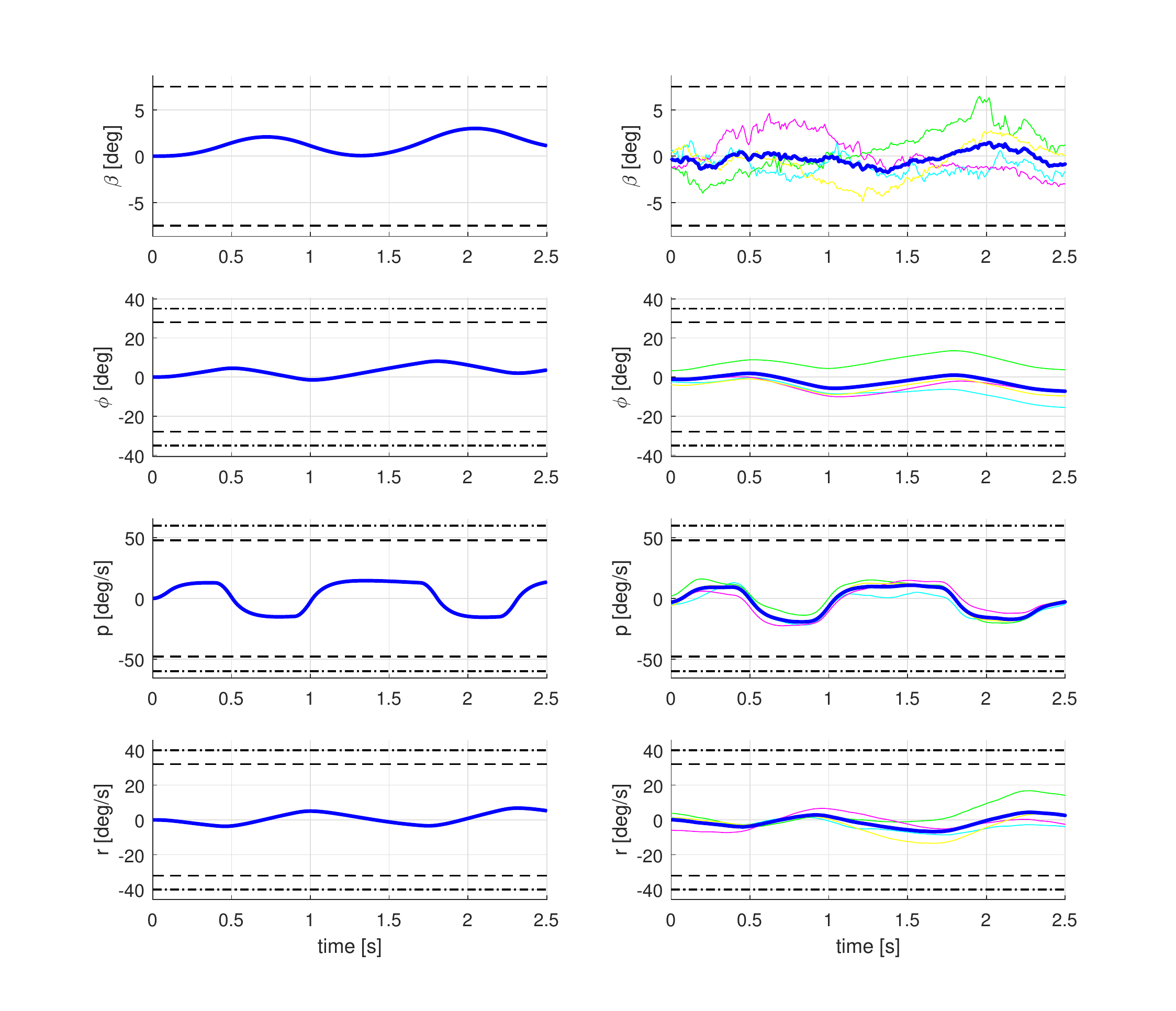}
	\caption{Simulated lateral response of the a priori model using optimal aileron deflection (left column) versus three real experiments (right column).  The real experiments are in green, magenta and cyan thin solid lines whereas their average is in thick blue solid line. In dash line the OED constraints while dash-dot line the flight envelope.}
	\label{fig:LatDyn_simVSreal_da}
\end{figure}

\begin{figure}[tbhp]
	\centering
	\includegraphics[width = 370pt, height = 300pt]{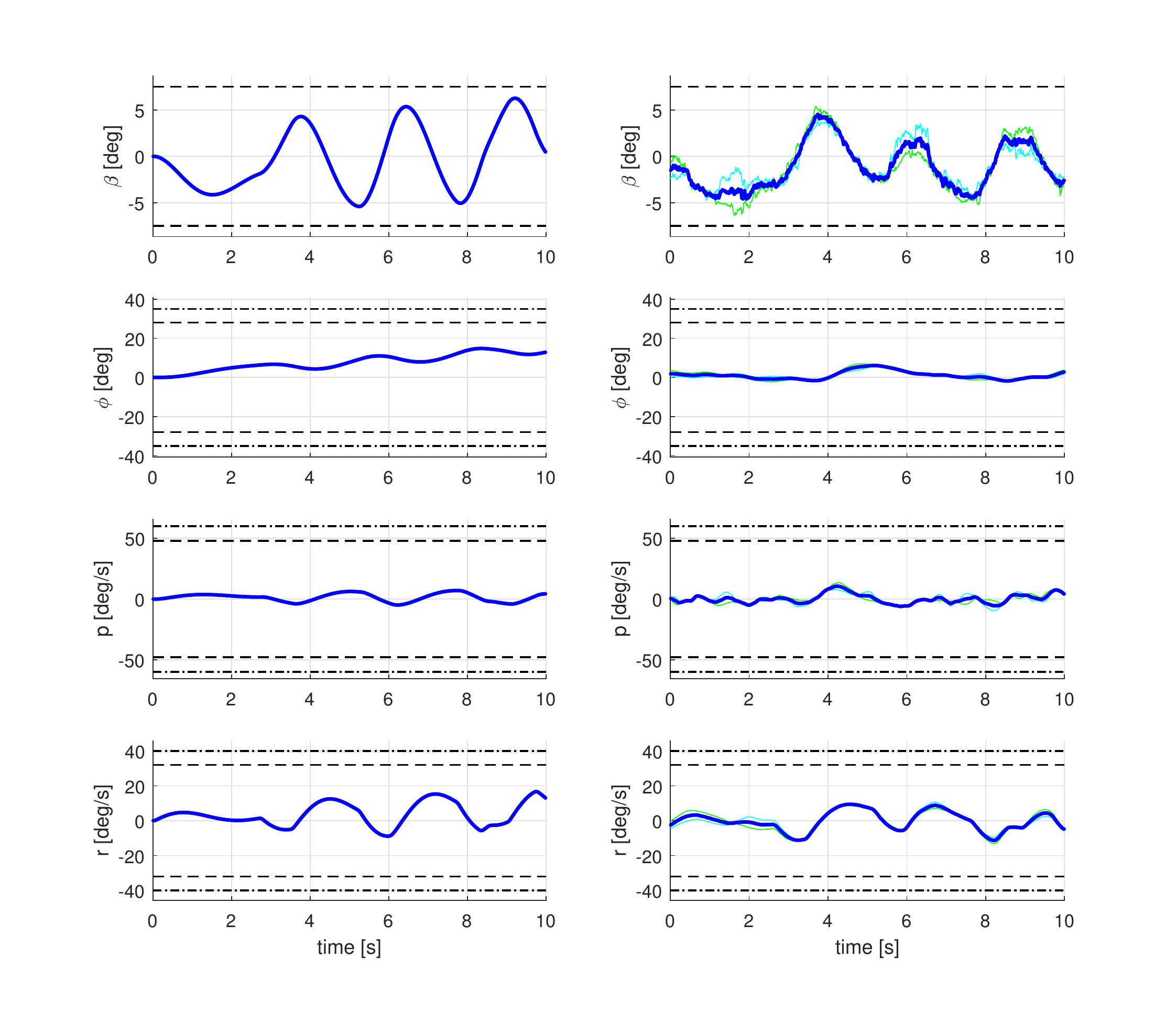}
	\caption{Simulated lateral response of the a priori model using optimal rudder deflection (left column) versus two real experiments (right column).  The real experiments are in green and cyan thin solid lines whereas their average is in thick blue solid line. In dash line the OED constraints while dash-dot line the flight envelope.}
	\label{fig:LatDyn_simVSreal_dr}
\end{figure}

The data obtained with the experiments are shown in figs. \ref{fig:LongDyn_simVSreal},\ref{fig:LatDyn_simVSreal_da},\ref{fig:LatDyn_simVSreal_dr} in comparison to the corresponding simulated optimal responses. The control surfaces deflection are omitted since they are already shown in section \ref{sec:SimulationResult}. A video that documents the system identification flight test is available at \cite{OptSysIDvideo}. Because of the difficulties to obtain accurate estimations of the a priori Phugoid mode, the airspeed response exceed the \ac{OED} constraints of $\approx 3 \, \mathrm{m/s}$ while the remaining responses in fig.~ \ref{fig:LongDyn_simVSreal} are bounded.

In fig.~\ref{fig:LatDyn_simVSreal_da}, it is shown that within two of three experiments the roll angle $\phi$ oscillates and moves towards the flight envelope limits whereas the yaw dynamics, i.\,e., $\beta$ and $r$, are barely excited. Opposite situations occur when the optimal sequence on the rudder control surface is performed, see fig.~\ref{fig:LatDyn_simVSreal_dr}. In this case, the roll angle is maintained close to its trim position along the entire open-loop sequence, showing less cross-coupling than predicted. Finally, the sinusoidal motion of the side slip $\beta$ does not increase in terms of amplitude as observed in simulation. In both lateral excitations the constraints were thoroughly fulfilled.

In summary, experimental data have shown that optimized maneuvers obtained with the proposed method are able to significantly excite the aircraft dynamics for a relatively long time window without violation of the flight envelope.  

\section{Conclusions}\label{sec:Conclusions}
In this paper, a time domain, model-based \acf{OED} and subsequent flight tests have been carried out for an autonomous aircraft. A suitable and comprehensive non-linear mathematical model for OED problems related to aircraft dynamics was introduced and an overview of the flight test procedure has been provided. The optimized maneuvers were obtained separately for the roll, pitch and yaw dynamics for the steady state wing level trim condition. The optimization problem was initialized using a priori aerodynamic derivatives obtained via lifting line method augmenting them with previous flight test campaign. The optimal solutions were compared with the well known and widely used 3-2-1-1 maneuvers and their estimation performance were assessed by the Cramen Rao Lower Bound. Simulation results have shown that optimized maneuvers enhance the information content of the experiment and help to validate and refine the a priori models. Finally, real flight test have been successfully conducted showing that given a fair accuracy of the a priori model, the open-loop aircraft response is bounded in the prescribed constraints. However, for safety reasons the optimal aileron sequence was not carried out completely.  


\section{Future works}\label{sec:FutureWorks}
Future works will aim to the investigation of feasible optimal aileron sequences which can safely excite the aircraft roll axis. 

\section*{Acknowledgments}
The authors would like to thank Ampyx Power B.V., and in particular the flight operations team for conducting the system identification flight tests.

\section*{Tables}
\begin{table}[tbhp]
	\caption{A priori longitudinal dimensional aerodynamic derivatives}
	\label{tab:LongitudinalDerivatives}
	\centering
	\begin{tabular}{ cc cc cc }\hline
		\bf X-axis       &\bf Value  & \bf Z-axis              &\bf Value &\bf M-axis        &\bf Value \\\hline
		$X_{V}$          &  -0.147   & $Z_{V}$                 & -0.060   & $M_{V}$          &   0.0    \\
		$X_{\alpha}$     &   7.920   & $Z_{\alpha} / \VTe$     & -4.400   & $M_{\alpha}$     &  -6.180  \\
		$X_{q}$          &  -0.163   & $Z_{q}$                 &  0.896   & $M_{q}$          &  -1.767  \\
		$X_{\delta_{e}}$ &  -0.232   & $Z_{\delta_{e}} / \VTe$ & -0.283   & $M_{\delta_{e}}$ & -10.668  \\\hline		
	\end{tabular}
\end{table}
\begin{table}[tbhp]
	\caption{A priori lateral dimensional aerodynamic derivatives}
	\label{tab:LateralDerivatives}
	\centering
	\begin{tabular}{ cc cc cc }\hline
		\bf Y-axis             &\bf Value &\bf L-axis         &\bf Value &\bf N-axis         &\bf Value \\\hline
		$Y_{\beta}/\VTe$       & -0.167   & $L_{\beta}'$      &  -8.201  & $N_{\beta}'$      &  3.214   \\
		$Y_{p}   $             &  0.0     & $L_{p}'$          & -11.292  & $N_{p}'$          & -0.750   \\
		$Y_{r}   $             & -0.976   & $L_{r}'$          &   3.853  & $N_{r}'$          & -0.457   \\
		$Y_{\delta_{a}}/ \VTe$ & -0.046   & $L_{\delta_{a}}'$ & -32.600  & $N_{\delta_{a}}'$ &  0.716   \\
		$Y_{\delta_{r}}/ \VTe$ &  0.093   & $L_{\delta_{r}}'$ &   0.524  & $N_{\delta_{r}}'$ & -2.370   \\\hline		
	\end{tabular}
\end{table}
\begin{table}[tbhp]
	\caption{Sensors noise standard deviation $\sigma_{y}$}
	\label{tab:std_sensors}
	\centering
	\begin{tabular}{lccc}\hline
		\bf Sensor           & \bf Variable         & $\mathbf{\sigma_{y}}$ & \bf Unit             \\ \hline
		five hole pitot tube & $V_{T}$              &    $2.5$     & $[\mathrm{m/s}]   $ \\
		five hole pitot tube & $(\alpha,\beta)$     &    $0.5$     & $[\mathrm{deg}]   $ \\
		IMU	                 & $(\phi,\theta,\psi)$ &    $0.1$     & $[\mathrm{deg}]   $ \\ 
		IMU	                 & $(p,q,r)$            &    $0.1$     & $[\mathrm{deg/s}] $ \\\hline		
	\end{tabular}
\end{table}

\bibliography{mybibfile}

\end{document}